\documentclass[11pt]{elsarticle}
\usepackage{fancyhdr}
\usepackage{marginnote}

\usepackage{amsmath}
\usepackage{mathtools}
\usepackage{amsfonts}
\usepackage{amsthm}
\usepackage{amssymb}
\usepackage{color}
\usepackage{dsfont}
\usepackage{geometry}
\newgeometry{tmargin=2.8cm, bmargin=3.5cm, lmargin=1.7cm, rmargin=1.7cm}

\usepackage{color}

\definecolor{darkgreen}{rgb}{0.00, 0.50, 0.00}

 \def\thebibliography#1{\section*{References\markboth
  {References}{References}}\list
  {[\arabic{enumi}]}{\settowidth\labelwidth{[#1]}\leftmargin\labelwidth
  \advance\leftmargin\labelsep
  \usecounter{enumi}}
  \def\newblock{\hskip .11em plus .33em minus -.07em}
  \sloppy
  \sfcode`\.=1000\relax}

 \let\OLDthebibliography\thebibliography
 \renewcommand\thebibliography[1]{
   \OLDthebibliography{#1}
   \setlength{\parskip}{3pt}
   \setlength{\itemsep}{0pt plus 1\baselineskip}
 }

\newtheorem{theo}{\bf Theorem} 
 
\newtheorem{coro}{\bf Corollary}[section]
\newtheorem{lem}[coro]{\bf Lemma} 
\newtheorem{rem}[coro]{\bf Remark}

\newtheorem{prop}[coro]{\bf Proposition}

\def\R{{\mathbb{R}}}

\def\r{{\mathbb{R}}}

\def\rn{{\mathbb{R}^{n}}}
\def\rd{{\mathbb{R}^{d}}}
\def\rnd{{\mathbb{R}^{n \times d}}}
\def\Rn{{\rn}}

\def\rp{{[0,\infty )}}

\def\ve{{\varepsilon}}
\def\vr{{\varrho}}

\def\MBm{{M^{-}_{B}}}
\def\MBxd{{\left(M^{-}_{B_{\tau\delta}(x)}\right)^{**}}}
\def\MBp{{M^{+}_{B}}}
\def\Msm{{(\MBm)^{**}}}

\def\Bcong{{\mathsf{(B^{gen}_\gamma)}}}
\def\Bcongp{{\mathsf{(B^{gen}_{\gamma+})}}}
\def\Bcongi{{\mathsf{(B^{iso}_\gamma)}}}
\def\Bcongo{{\mathsf{(B^{ort}_\gamma)}}}

\def\bu{{\mathbf{u}}}
\def\bxi{{\boldsymbol{\xi}}}
\def\bvp{{\boldsymbol{\vp}}}
\def\bvpd{{\boldsymbol{\vp_\delta}}}

\def\iO{{\int_{\Omega}}}

\newcommand{\vp}{\varphi}

\newcommand{\dv}{\mathrm{div}}
\newcommand{\cF}{{\mathcal{F}}}
\newcommand{\cG}{{\mathcal{G}}}
\newcommand{\cH}{{\mathcal{H}}}

\begin{document}

\begin{frontmatter}

\title{Absence of Lavrentiev's gap for anisotropic functionals}

\author[1]{Michał Borowski}\ead{m.borowski@mimuw.edu.pl}
\author[1]{Iwona Chlebicka\corref{mycorrespondingauthor}}
\cortext[mycorrespondingauthor]{Corresponding author}
\ead{i.chlebicka@mimuw.edu.pl} 
\author[1]{Błażej Miasojedow}\ead{b.miasojedow@mimuw.edu.pl}

\fntext[myfootnote]{MSC2020: 46E30 (46E40,46A80)}
\fntext[myfootnote]{M.B. and I.C. are supported by NCN grant 2019/34/E/ST1/00120.}

\address[1]{Institute of Applied Mathematics and Mechanics,
University of Warsaw, ul. Banacha 2, 02-097 Warsaw, Poland
}

\begin{abstract} 
We establish the absence of the Lavrentiev gap between Sobolev and smooth maps for a non-autonomous variational problem of a general structure, where the integrand is assumed to be controlled by a function which is convex and anisotropic with respect to the last variable. This fact results from new results on good approximation properties of the natural underlying unconventional function space. {Scalar and vector-valued problems are studied.} 
\end{abstract}

\end{frontmatter}


\section{Introduction}
We study the well-posedness of minimization of the following variational functional
\begin{align}\label{cF-def}
    \cF[u]:=\int_\Omega F(x,u,\nabla u)\,dx\,,
\end{align}
over an open and bounded set $\Omega\subset\rn$, $n\geq 1$, where $F:\Omega\times\R\times\rn\to\R$ is merely continuous with respect to the second and the third variable. We suppose that there exist constants $0<\nu,\beta<1<L$ such that
\begin{align}\label{sandwich}
    \nu M\left(x,{\beta} \xi\right)\leq F(x,z,\xi)\leq L( M\left(x,\xi\right)+1), \qquad\text{for all }\ x\in\Omega,\ z\in\R,\ \xi\in\rn, 
\end{align}
for an $N$-function $M:\Omega\times\rn\to\rp$ that is continuous with respect to both variables. This function is called {\it anisotropic} because it depends on $\xi$ not necessarily via $|\xi|$. We focus on the issue whether the minimizers can be approximated by regular functions in the topology that is natural for the problem. As an answer, we find sharp conditions on $M$ ensuring the absence of the Lavrentiev gap between Sobolev and smooth maps that precisely capture the anisotropic version of double-phase and multi-phase functionals. This fact results from fine approximation properties of the natural underlying unconventional function space.

The natural function space to study minimizers to problems like~\eqref{cF-def} is a Sobolev type space equipped with a Luxemburg norm defined by the means of a functional $\xi\mapsto\int_\Omega M(x,\xi)\,dx$ applied to a distributional gradient of a function $u\in W^{1,1}(\Omega)$ with $\nabla u$ in a relevant Musielak--Orlicz spaces. Since the growth of $M$ does not need to  be doubling, it should be taken into account that the Musielak--Orlicz spaces $E_M$ and $L_M$ (being norm and modular closure of $L^\infty$, respectively)  may differ. Consequently, their Sobolev-type versions $V^1_0 E_M(\Omega)$ and $V^1_0 L_M(\Omega)$ do not need to coincide. If $M,M^*\in\Delta_2$, then $L_M=E_M$ and $V^1_0 L_M=V^1_0 E_M=V_0^{1,M}$. 

 The situation when the infimum of a variational problem over a family of regular functions (e.g. smooth or Lipschitz) is strictly greater than the infimum taken over all functions satisfying the same boundary conditions 
 is called Lavrentiev's phenomenon after his seminal paper~\cite{Lav}. To give a flavour of a great deal of classical results and new developments in this area, let us refer to~\cite{badi,basu,Bor-Chl,Bousquet-Pisa,bgs-arma,Buttazzo-Belloni,Buttazzo-Mizel,DeFilippisLeonetti+2022,eslev,eslemi,MiRa,zh86} and references therein. In our setting,  having $u_0\in V^1E_M(\Omega)$, if $M$ is not sufficiently regular, it might happen that
\begin{align}
    \inf_{u_0+V^1_0E_M(\Omega)}\cF[u]<\inf_{u_0+C_c^\infty(\Omega)}\cF[u]\,.\label{gap}
\end{align}
Since~\cite{zh86, zh-dens} by Zhikov, it is known that  log-H\"older continuity of the variable exponent of $M(x,\xi)=|\xi|^{p(x)}$ ensures that the above scenario is excluded. In other kinds of isotropic generalized Orlicz growth problems, the role of this continuity condition is typically played by an assumption~(A1') from~\cite{hahab,ha-jfa}, embracing conditions from~\cite{bacomi-cv, comi,DeFilippisLeonetti+2022, eslev} needed for study on  $M(x,\xi)=|\xi|^{p}+a(x)|\xi|^q$ with possibly vanishing weight $a$. Such conditions are inevitable due to the examples of functions that cannot be approximated~\cite{badi,eslemi,fomami,Marc-ARMA-89,zh}. Smooth approximation properties of an inhomogeneous and anisotropic space of Musielak--Orlicz--Sobolev-type were proven under some far from sharp conditions~\cite{C-b,cgzg,gszg} and improved recently to a truly local and anisotropic one in~\cite{Bor-Chl}. Nonetheless, in the view of~\cite{badi,bgs,bgs-arma}, there was still some room for progress. Apart from the intrinsic mathematical interest of Lavrentiev's phenomenon, we shall indicate that the Musielak--Orlicz--Sobolev spaces are used as a framework in modelling of modern materials involving electrorheological and non-Newtonian fluids, thermo-visco-elastic ones, as well as in image restoration processing, see~\cite{IC-pocket,PA,klawe,AW}.  Henceforth, we provide a precise tool into their analysis. Section~\ref{sec:appl-PDEs} describes how our result directly extend the existence results of~\cite{C-b,gszg} and other contributions.\newline

One of the most important feature of the space that we want to include in our analysis is anisotropy. A weak $N$-function $M:\Omega\times\rn\to[0,\infty)$ is called {\it isotropic} if $M(x,\xi)=m(x,|\xi|)$ with a weak $N$-function $m:\Omega\times\rp\to\rp$ and {\it anisotropic} if its dependence on $\xi$ is allowed to be more complicated. One can consider the functions admitting a decomposition called {\it orthotropic} (studied e.g. in~\cite{BB,DFG}): \[\text{$M\big(x,(\xi_1,\dots,\xi_n)\big)=\sum_{i=1}^n M_i(x,|\xi_i|)\quad$ with weak $N$-functions $\ M_i:\rp\to\rp$.}\] If an anisotropic function is not comparable to any function admitting a decomposition that after an affine and invertible change of variables has the above form, we call it {\it essentially fully anisotropic}. A relevant example can be found in~\cite{cn}.  Fully anisotropic spaces  are considered since~\cite{Ioffe,Klimov76,Sk1} and~\cite{BGM,Cfully,Sch}. They serve as a setting for nonlinear partial differential equations or calculus of variations, see e.g.~\cite{ACCZG,barlettacianchi,Ci-sym,Stroff}. Recently, we can observe more and more attention paid to problems that are both inhomogeneous and anisotropic at the same time, e.g.~\cite{bgs,C-b,ren-para,cgzg,CKL,gszg,PA,ha-aniso,Li,AW} and sharp conditions for the density of smooth functions are strikingly missing in the theory.\newline

Our goal is to detect the optimal conditions that need to be imposed on $M$ to exclude the case of strict inequality in~\eqref{gap}. The key accomplishments of the present paper yield the absence of Lavrentiev's phenomenon (Theorem~\ref{theo:lavrentiev}) and  the modular density of smooth functions (Theorem~\ref{theo:approx}). Both of the mentioned results holds for weak $N$-functions satisfying some balance condition reflecting a log-H\"older continuity of the variable exponent of a sharp range of powers in the double-phase case. In order to present the conditions, given a weak $N$-function $M:\Omega\times \rn\to\rp$ and a ball $B$, we define
\begin{align}\label{eq:defMpMmm}
    \MBm(\xi) = {\rm ess\,inf}_{y\in B}M(y, \xi)\quad\text{and}\quad
    \MBp(\xi) = {\rm ess\,sup}_{y\in B}M(y,\xi)\,.
\end{align}
 The considered conditions read as follows. 

\begin{description}
    
\item {\bf Isotropic condition $\Bcongi$.} Let us assume that $M:\Omega\times\rp\to\rp$ is a weak $N$-function. Suppose there exist constants $c^\Diamond,C^\Diamond\geq1$ such that 
for every ball $B$ with radius $r \leq 1$ and for all $\xi\in\rn$ satisfying $|\xi|\leq c^\Diamond r^{\gamma-1}$ there holds $M_{B}^+(|\xi|) \leq M_{B}^-(C^\Diamond|\xi|) + 1$.

\item {\bf Orthotropic condition $\Bcongo$.} Let us assume that 
\[
M(x,\xi)=\sum_{i=1}^n M_i(x,|\xi_{i}|),
\]
where $M_i:\Omega\times\rp\to\rp$ are weak $N$-functions and $\xi=(\xi_1,\dots,\xi_n)$. Suppose there exist constants $c^\Diamond,C^\Diamond\geq1$ such that 
for every ball $B$ with radius $r \leq 1$, all $i$, and for all $\xi_i\in\R$ satisfying $|\xi_i|\leq c^\Diamond r^{\gamma-1}$  there holds $(M_i)_{B}^+(|\xi_i|) \leq (M_i)_{B}^-(C^\Diamond|\xi_i|) + 1$.
\end{description}

\noindent The above conditions $\Bcongi$ and $\Bcongo$ covering the most typical settings are special cases of a fully anisotropic condition $\Bcong$, see Section~\ref{sec:general-growth}. They  force that the growth of a  weak $N$-function with respect to the second variable is not perturbed much in a small spacial region. The first of our main accomplishments yields precise result on the absence of Lavrentiev's phenomenon, i.e.,
\begin{align}
    \inf_{u_0+V^1_0L_M(\Omega)}\cF[u]=\inf_{u_0+C_c^\infty(\Omega)}\cF[u]\,,\label{no-gap-intro}
\end{align} whereas the second one -- the approximation result in anisotropic Musielak--Orlicz--Sobolev spaces, i.e.,
\begin{align*}
V_0^1L_M(\Omega)=\overline{C_c^\infty(\Omega)}^{mod}
\,.
\end{align*}Here $\overline{C_c^\infty(\Omega)}^{mod}$ stands for the closure in the sequential modular topology of $L_M$  of the gradients. Theorem~\ref{theo:lavrentiev} is supplied with its extended versions in Section~\ref{sec:general-growth}, which are more complicated in the exposition, but cover more general functionals. We embrace and extend the results on the absence of Lavrentiev's phenomenon from~\cite{yags,bgs-arma,eslemi} to precisely capture local nature of the problem, as well as anisotropy and general growth of $F$ with respect to the last variable. Moreover, in the regions of low growth of $M$ we improve the known anisotropic approximation results of \cite{Bor-Chl,C-b,gszg}. When the growth is isotropic, we do it up to the known borderline cases, see~\cite{basu,bacomi-cv,fomami,zh86}. In turn, we supply the existence theory that holds in the absence of Lavrentiev's phenomenon~\cite{bgs,C-b,ren-para,cgzg,CKL,gszg} with precise information when the methods therein apply. Let us show two simple examples directly illustrating how our main results extend the state of the art. \newline

\noindent {\bf Main models.} A consequence of our main result, that is Theorem~\ref{theo:lavrentiev}, is that for a double-phase functional \begin{align}\label{def-double-phase-functional}
    \cH_{\rm iso}[u]:=\int_\Omega b(x,u)\left(|\nabla u|^p+a(x) |\nabla u|^{q}\right)\,dx
\end{align} where $b:\Omega\times\r\to\r$ satisfies $0<\nu<b(\cdot,\cdot)<L$ and is merely continuous with respect to the last variable, $1\leq p\leq q$, $a:\Omega\to[0,\infty)$ such that $a\in C^{0,\alpha}(\Omega)$, and for
\[H_{\rm iso}(x,\xi)=|\xi|^p+a(x) |\xi|^{q},\quad \text{and}\quad V^{1,H_{\rm iso}}_0(\Omega):=\{u\in W^{1,1}_0(\Omega):\ |\nabla u|\in L^{H_{\rm iso}}(\Omega)\},\] it holds
\begin{align}
    \inf_{u_0+V_0^{1,H_{\rm iso}}(\Omega)}\cH_{\rm iso}[u]=\inf_{u_0+C_c^\infty(\Omega)}\cH_{\rm iso}[u]\,,\label{no-gap-dp}
\end{align}
whenever\begin{align}
    q\leq p+\alpha\,.\label{qleqpplusalpha}
\end{align}
Consequently, even in the classical isotropic case, Theorem~\ref{theo:lavrentiev} extends the known results. In particular, it generalizes the very recent contribution~\cite{bgs-arma} relaxing growth by allowing for functionals involving integrands dependent on three variables (precisely, $\cF$ like~\eqref{cF-def} under the assumption~\eqref{sandwich} with $M$ substituted by $H_{\rm iso}$). In the appearance of more phases, we give more precise bound than~\cite{bgs-arma}. Furthermore, having a priori knowledge on the regularity of a minimizer Theorem~\ref{theo:lavrentiev} enables to improve the range~\eqref{qleqpplusalpha}. Namely, for $u,u_0\in C^{0,\gamma}(\Omega)$, $\gamma\in[0,1]$, we get that~\eqref{no-gap-dp} holds whenever
\[q\leq p+\frac{\alpha}{1-\gamma}\,.\]
Since we need only continuity of $b$ and only with respect to the second variable, we relax the assumptions of \cite[Theorem~4]{bacomi-cv} yielding the absence of Lavrentiev's gap for $\cH_{\rm iso}$ under extra assumptions on the decay of the modulus of continuity of $b$ with respect to both variables. The range from~\eqref{qleqpplusalpha} cannot be improved for $p<n$ due to~\cite{badi}. Our main model, that was not covered by the literature, is the following anisotropic functional
\begin{align}
    \cG[u]:=\int_\Omega G(x,u,\nabla u)\,dx\,,
\end{align}where $G:\Omega\times\R\times\rn\to\R$ is continuous with respect to all its variables and \begin{align*}
\nu H\left(x,\xi\right)\leq    G(x,z,\xi)\leq L H\left(x,\xi\right),\quad\text{for all}\ x\in\Omega,\ z\in\R,\ \xi\in\rn\,
\end{align*}
and some constants $\nu, L>0$, where
\begin{align}\label{eq:def-mp}
    \text{ $H(x,\xi)=\sum_{i=1}^n |\xi_i|^{p_i}+\sum_{i=1}^n a_i(x) |\xi_i |^{q_i},$ $1{\leq} p_i\leq q_i$, $a_i:\Omega\to[0,\infty)$,  such that $a_i\in C^{0,\alpha_i}(\Omega)$},\ i=1,\dots,n\,.
\end{align}
As a consequence of our main result, it holds
\begin{align}\label{aniso-double-phase-no-gap}
    \inf_{u_0+V_0^{1,H}(\Omega)}\cG[u]=\inf_{u_0+C_c^\infty(\Omega)}\cG[u]\,
\end{align}
\[\text{where } V^{1,H}_0(\Omega):=\{u\in W^{1,1}_0(\Omega):\ |\nabla u|\in L^{H}(\Omega)\},\]
whenever\begin{align}\label{aniso-double-phase-range}
    q_i\leq p_i+\alpha_i \quad\text{for }\ i=1,\dots,n\,.
\end{align}
Furthermore, in this case we can also trade a priori known  regularity of a minimizer with this range. For $u,u_0\in C^{0,\gamma}(\Omega)$, $\gamma\in[0,1]$, by Theorem~\ref{theo:lavrentiev} we get that~\eqref{aniso-double-phase-no-gap} holds whenever
\[q_i\leq p_i+\frac{\alpha_i}{1-\gamma} \quad\text{for }\ i=1,\dots,n\,.\]
So far, the best known result for the absence of Lavrentiev's phenomenon in this kind of anisotropic double-phase spaces was due to~\cite{Bor-Chl} and covered the ranges for the exponents $\frac{q_i}{p_i}\leq 1+\frac{\alpha_i}{n}$, $i=1,\dots,n$. Note that when $p_i<n$, this range of admissible $p_i$ and $q_i$ is smaller than~\eqref{aniso-double-phase-range}. This means that in the current study, we admit a worse modulus of continuity of $H$ and still provide \eqref{aniso-double-phase-no-gap}. See also Remark~\ref{rem:pleqn}.

In Corollary~\ref{coro:lavrentiev}, we give more examples where the modulus of continuity is essentially improved compared to~\cite{Bor-Chl}, including anisotropic variable exponent double-phase functionals and functionals Orlicz phases. Remark~\ref{rem:gen-growth} illustrates our results in the setting of the general growth and full anisotropy. To our best knowledge, there is no anisotropic counterexample available so far in the literature.\newline


    Let us pass to presenting our main accomplishments.  Function space $L_M$ is defined in Section~\ref{sec:prelim}. Our main results concern \begin{align*}
    V_0^{1}L_M(\Omega)&=\big\{f\in W^{1,1}_{0}(\Omega): \ \nabla f\in L_M(\Omega;\rn)\big\}\,.
\end{align*}

\noindent{\bf Absence of Lavrentiev's phenomenon.}
Let us formulate our general result that under a balance condition Lavrentiev's phenomenon for $\cF$ does not occur between Sobolev and smooth maps. It is followed by a long list of examples being of separate attention in the field.
\begin{theo}[Absence of Lavrentiev's phenomenon]
\label{theo:lavrentiev}
Let $\Omega$ be a bounded Lipschitz domain in $\rn$ and functional $\cF$ be given by~\eqref{cF-def} with $F:\Omega\times\R\times\rn\to\R$ satisfying~\eqref{sandwich} for a weak $N$-function $M$ that is continuous with respect to both variables and such that $M\in\Delta_2$. Assume further that $F$ is measurable with respect to the first variable and continuous with respect to the second and the third variable. Then we observe the absence of Lavrentiev's phenomenon in the following cases.
\begin{enumerate}[{\it (i)}]
    \item If  $\gamma=0$, $u_0\in V^1L_M(\Omega)$, and $M$ satisfies condition $\Bcongi$ or $\Bcongo$, then\begin{align}
    \inf_{u_0+V^1_0L_M(\Omega)}\cF[u]=\inf_{u_0+C_c^\infty(\Omega)}\cF[u]\,.\label{no-gap-gen}
\end{align} 
    \item If $\gamma\in(0,1)$, $u_0\in V^1L_M(\Omega)\cap C^{0,\gamma}(\Omega)$, and $M$ satisfies condition $\Bcongi$ or $\Bcongo$, then
\begin{align}
    \inf_{u_0+V^1_0L_M(\Omega)\cap C^{0,\gamma}(\Omega)}\cF[u]=\inf_{u_0+C_c^\infty(\Omega)}\cF[u]\,.\label{no-gap-hold}
\end{align} 
    \item If $u_0\in V^1L_M(\Omega)\cap C^{0,1}(\Omega)$, then
\begin{align}
    \inf_{u_0+V^1_0L_M(\Omega)\cap C^{0,1}(\Omega)}\cF[u]=\inf_{u_0+C_c^\infty(\Omega)}\cF[u]\,.\label{no-gap-lip}
\end{align} 
\end{enumerate}
\end{theo}

\begin{rem}Conditions $\Bcongi$ and $\Bcongo$ are always satisfied when $M(x,\xi)=M(\xi)$ including anisotropic functions. 
In turn, functionals driven by such $M$ never face Lavrentiev's phenomenon.
\end{rem}

By a direct application of Theorem~\ref{theo:lavrentiev} in the case of particular choices of $M$ we get the following corollary.
\begin{coro}
\label{coro:lavrentiev} 
Let $\Omega$ be a bounded Lipschitz domain in $\rn$, $\gamma\in [0,1)$, and functional $\cF$ be given by~\eqref{cF-def} with $F:\Omega\times\R\times\rn\to\R$ satisfying~\eqref{sandwich} for a weak $N$-function $M$. Assume further that $F$ is measurable with respect to the first variable and continuous with respect to the second and the third variable. Then we have the following isotropic consequences of Theorem~\ref{theo:lavrentiev}.
\begin{enumerate}[{\it (i)}]
    \item (variable exponent) If $M(x,\xi)=|\xi|^{p(x)}$, $p\in {\cal P}^{\log}(\Omega)$, $1\leq p(\cdot)<\infty$, then Lavrentiev's phenomenon for $\cF$ does not occur between $W^{1,p(\cdot)}_0(\Omega)$ and $C_c^\infty(\Omega)$ {as well as between $W^{1,p(\cdot)}_0(\Omega)\cap C^{0,\gamma}(\Omega)$ and $C_c^\infty(\Omega)$}. See~\cite{CUF,DHHR,zh86}.   
    \item (mild double-phase) If $M(x,\xi)=|\xi|^{p}+a(x)|\xi|^p\log(e+|\xi|)$, $0\leq a\in {\cal P}^{\log}(\Omega)$, $1\leq p<\infty$, then Lavrentiev's phenomenon for $\cF$ does not occur between $V^{1,M}_0(\Omega)$ and $C_c^\infty(\Omega)$ {as well as between $W^{1,{M}}_0(\Omega)\cap C^{0,\gamma}(\Omega)$ and $C_c^\infty(\Omega)$}. See~\cite{bacomi-st}.
    \item (double-phase) If $M(x,\xi)=|\xi|^{p}+a(x)|\xi|^q$, $0\leq a\in C^{0,\alpha}(\Omega)$, $1\leq p\leq q$, then Lavrentiev's phenomenon for $\cF$ does not occur\begin{itemize}[--]
        \item between $V^{1,M}_0(\Omega)$ and $C_c^\infty(\Omega)$ whenever $q\leq p+\alpha$,
        \item between $V^{1,M}_0(\Omega)\cap C^{0,\gamma}(\Omega)$ and $C_c^\infty(\Omega)$ whenever $q\leq p+\alpha/(1-\gamma)$.
    \end{itemize}  See~\cite{bacomi-cv,comi, DeFilippisLeonetti+2022,eslemi,eslev}.
    \item (variable exponent double-phase) If $M(x,\xi)=|\xi|^{p(x)}+a(x)|\xi|^{q(x)}$, $0\leq a\in C^{0,\alpha}(\Omega)$, $1\leq p(x)\leq q(x)$, {$p, q \in \mathcal{P}^{log}(\Omega)$} then Lavrentiev's phenomenon for $\cF$ does not occur \begin{itemize}[--]
        \item between $V^{1,M}_0(\Omega)$ and $C_c^\infty(\Omega)$ whenever $\sup_{x\in\Omega} (q(x)-p(x))\leq\alpha$,
        \item between $V^{1,M}_0(\Omega)\cap C^{0,\gamma}(\Omega)$ and $C_c^\infty(\Omega)$ whenever  $\sup_{x\in\Omega} (q(x)-p(x))\leq \alpha/(1-\gamma)$.
    \end{itemize}  See~\cite{bgs-arma,rata}. 
    \item (multi-phase) If $M(x,\xi)=|\xi|^{p}+\sum_{i=1}^n a_i(x)|\xi|^{q_i}$, $0\leq a_i\in C^{0,\alpha_i}(\Omega)$, $1\leq p\leq q_i$, $i=1,\dots,n$, then Lavrentiev's phenomenon for $\cF$ does not occur \begin{itemize}[--]
        \item between $V^{1,M}_0(\Omega)$ and $C_c^\infty(\Omega)$ whenever  $q_i\leq p+\alpha_i$, 
        \item between $V^{1,M}_0(\Omega)\cap C^{0,\gamma}(\Omega)$ and $C_c^\infty(\Omega)$ whenever  $q_i\leq p+\alpha_i/(1-\gamma)$.  
    \end{itemize}  See~\cite{DeF-multi}.
    \item (Orlicz double-phase) If $M(x,\xi)=\phi(|\xi|)+ a(x)\psi(|\xi|)$, $\phi,\psi\in\Delta_2$, $\psi(t)/\phi(t)\xrightarrow[t\to\infty]{}\infty$, $0\leq a\in C(\Omega)$ has a modulus of continuity $\omega_a$, then Lavrentiev's phenomenon for $\cF$ does not occur  
    \begin{itemize}[--]
        \item between $V^{1,M}_0(\Omega)$ and $C_c^\infty(\Omega)$ whenever  $\omega_a$ satisfies {$\omega_{a}(t) \leq \tfrac{\phi(t^{-1})}{\psi(t^{-1})}$},
        \item between $V^{1,M}_0(\Omega)\cap C^{0,\gamma}(\Omega)$ and $C_c^\infty(\Omega)$ whenever  $\omega_a$ satisfies {$\omega_{a}(t) \leq \tfrac{\phi(t^{\gamma-1})}{\psi(t^{\gamma-1})}$}.
    \end{itemize}
    In particular, if $\phi(t) = t^p\log^{-\beta}(e + t)$, $\psi(t) = t^p\log^\alpha(e + t)$, $1<p<\infty$, $\alpha,\beta\geq 0$, then the choice of $a$ satisfying $\omega_a(t) \leq {\log^{-(\alpha + \beta)}(e + t)}$ excludes the Lavrentiev's phenomenon, cf.\cite[Corollary 4]{basu}.\\
    For a refinement of this result without assumption $\phi,\psi\in\Delta_2$, see Section~\ref{sec:general-growth}.
\end{enumerate}
Moreover, we have the following anisotropic consequences of Theorem~\ref{theo:lavrentiev}.
    \begin{enumerate}[{\it (i)}]
    \item If $M(x,\xi)=\sum_{i=1}^n|\xi_i|^{p_i(x)}$, $p_i\in {\cal P}^{\log}(\Omega)$, $1\leq p_i(\cdot)<\infty$, $i=1,\dots,n$, $\vec{p}(x)=(p_1(x),\dots,p_n(x))$, then Lavrentiev's phenomenon for $\cF$ does not occur between $W^{1,\vec p(\cdot)}_0(\Omega)$ and $C_c^\infty(\Omega)$ {as well as between $W^{1,\vec p(\cdot)}_0(\Omega)\cap C^{0,\gamma}(\Omega)$ and $C_c^\infty(\Omega)$}.
    \item If $M(x,\xi)=\sum_{i=1}^n|\xi_i|^{p_i}+\sum_{i=1}^n a_i(x)|\xi_i|^{q_i}$, $0\leq a_i\in C^{0,\alpha_i}(\Omega)$, $1\leq p_i\leq q_i$, $q_i\leq p_i+\alpha_i$, $i=1,\dots,n$, then Lavrentiev's phenomenon for $\cF$ does not occur
    {\begin{itemize}[--]
        \item between $V^{1,M}_0(\Omega)$ and $C_c^\infty(\Omega)$ for $q_i\leq p_i+\alpha_i$,
        \item between $V^{1,M}_0(\Omega)\cap C^{0,\gamma}(\Omega)$ and $C_c^\infty(\Omega)$ for $q_i\leq p_i+\alpha_i/(1-\gamma)$.
    \end{itemize}}
    \item If $M(x,\xi)=\sum_{i=1}^n|\xi_i|^{p_i}+\sum_{i=1}^na_i(x)|\xi_i|^{p_i}\log(e+|\xi_i|)$, $0\leq a_i\in {\cal P}^{\log}(\Omega)$, $1\leq p_i<\infty$, $i=1,\dots,n$, then Lavrentiev's phenomenon for $\cF$ does not occur between $V^{1,M}_0(\Omega)$ and $C_c^\infty(\Omega)$ {as well as between $V^{1,M}_0(\Omega)\cap C^{0,\gamma}(\Omega)$ and $C_c^\infty(\Omega)$}.
    \item If $M(x,\xi)=\sum_{i=1}^n|\xi_i|^{p_i(x)}+\sum_{i=1}^na_i(x)|\xi_i|^{q_i(x)}$, $0\leq a_i\in C^{0,\alpha_i}(\Omega)$, $1\leq p_i(\cdot)\leq q_i(\cdot)$, {$p_i, q_i \in \mathcal{P}^{log}(\Omega)$,} $i=1,\dots,n$, then Lavrentiev's phenomenon for $\cF$ does not occur
    {\begin{itemize}[--]
        \item between $V^{1,M}_0(\Omega)$ and $C_c^\infty(\Omega)$ whenever $\sup_{x\in\Omega}(q_i(x)-p_i(x))\leq \alpha_i$,
        \item between $V^{1,M}_0(\Omega)\cap C^{0,\gamma}(\Omega)$ and $C_c^\infty(\Omega)$ whenever $\sup_{x\in\Omega}(q_i(x)-p_i(x))\leq\alpha_i/(1-\gamma)$.
    \end{itemize}}
\end{enumerate}
\end{coro}


\noindent{\bf Density results.} We will distinguish between two kinds of density  -- the modular density and the norm one. In order to do so, let us define the following two kinds of convergence. \newline

\noindent {\it We say that a sequence $\{\xi_k\}_{k=1}^\infty$ converges modularly to $\xi$ in~$L_M(\Omega;\rn)$ (and denote it by $\xi_k\xrightarrow[k\to\infty]{M}\xi$), if there exists $\lambda>0$ such that
\begin{equation*}
\int_{\Omega}M\left(x,\frac{\xi_k-\xi}{\lambda}\right)\,dx\to 0.
\end{equation*}We say that a sequence $\{\xi_k\}_{k=1}^\infty$ converges to $\xi$ in norm topology of~$L_M(\Omega;\rn)$, if $\|\xi_k-\xi\|_{L_M}\xrightarrow[k\to\infty]{} 0$.}\newline
 
The  topologies generated by modular and norm convergences coincide when $M$ grows regularly enough, that is $M\in\Delta_2$, but in general we only know that if $\|\xi_k-\xi\|_{L_M}\to 0$, then $\xi_k\xrightarrow[k\to\infty]{M}\xi$ in $L_M$.\newline

Our main result on approximation reads as follows.
\begin{theo}[Density of smooth functions]\label{theo:approx}
Let $\Omega$ be a bounded Lipschitz domain in $\rn$ and  $M$ be a weak $N$-function that is continuous with respect to both variables. Then the following assertions hold true. \begin{enumerate}[{\it (i)}]
\item If $\gamma=0$ and $M$ satisfies condition $\Bcongi$ or $\Bcongo$, then for any $\vp\in V_0^1L_M(\Omega)$  there exists a sequence $\{\vp_\delta\}_{\delta>0}\subset C_c^\infty(\Omega)$, such that  $ \vp_\delta\to  \vp$ strongly in $L^1(\Omega)$ and in measure, and $\nabla\vp_\delta\xrightarrow[\delta\to 0]{M}\nabla \vp$ modularly in $L_M(\Omega;\rn)$.
\item If  $\gamma\in(0,1]$ and $M$ satisfies condition $\Bcongi$ or $\Bcongo$, then  for any $\vp\in V_0^1L_M(\Omega)\cap C^{0,\gamma}(\Omega)$  there exists a sequence $\{\vp_\delta\}_{\delta>0}\subset C_c^\infty(\Omega)$, such that  $ \vp_\delta\to  \vp$ strongly in $L^1(\Omega)$ and in measure, and $\nabla\vp_\delta\xrightarrow[\delta\to 0]{M}\nabla \vp$ modularly in $L_M(\Omega;\rn)$.
\end{enumerate} 
Moreover, in both above cases, if $\vp\in L^\infty(\Omega)$, then there exists $c=c(\Omega)>0$, such that $\|\vp_\delta\|_{L^\infty(\Omega)}\leq c \|\vp \|_{L^\infty(\Omega)}$.
\end{theo}

\noindent{\bf Methods of the proofs.} The main tool of our proof is to employ convolution-based approximation, which is classical and widely used in many contexts. In particular, it is used in the proofs by Gossez~\cite{Gossez} of the smooth density in the modular topology in the classical Orlicz spaces and in the anisotropic version of this result~\cite{ACCZG}. We cannot apply it directly due to $x$-dependence of the weak $N$-function defining the underlying function space and the relevant (modular) convergence. In such situation one is forced to use some continuity properties of $M$, cf. e.g.~\cite{yags,eslemi,zh-dens}. The proof gets even more complicated in the presence of anisotropy. The search for optimal conditions resulted in many contributions, including~\cite{yags,Bor-Chl,bgs,C-b,gszg}. Typically, these attempts make use of the greatest convex minorant of local infima of $M$ and balancing its behaviour with some quantities that are possible to control. In the isotropic case, $M_B^-$ is `almost convex', that is, Jensen's inequality holds for this function with an intrinsic constant, see \cite[Lemma~4.3]{ha-jfa}. For anisotropic functions, such constant does not exist  in general. In fact, even if $x\mapsto M(x,\cdot)$ is regular, $(M_B^-)^{**}$ can be arbitrarily far from $M_B^-$. Due to the very recent results, it is possible to restrict the set on which good properties of $(M^-_B)^{**}$ are needed (see~\cite{Bor-Chl}) and to find a concise and verifiable condition on $M$ (see~\cite{ha-aniso}) that was sufficient for construction of approximation. A different approach to this part can be found in~\cite{bgs,bgs-arma}, where in the middle of a convolution-based approximation procedure the Young inequality for convolution is applied combined with the fact that differentiation commutes with convolution. This enabled to relax the necessary condition, in particular to capture the whole range $q\leq p+\alpha$ in double-phase functional from~\eqref{def-double-phase-functional} for $b\equiv 1$. Nonetheless, the isotropic and doubling approach of~\cite{bgs-arma} is not crafted to track the growth of $M$ locally, nor allows studying $b:\Omega\times\R\to\R$. Despite the crucial idea comes from anisotropic attempt of~\cite{bgs}, the condition therein does not keep precise control on anisotropy. 

We extend applicability of all the mentioned contributions, i.e., \cite{yags,Bor-Chl,bgs,bgs-arma,C-b}, to the borderline cases and to capture phenomena that have local, general growth, and anisotropic nature in one shot. The main novelty of our approach, that allows to relax all previous conditions and additionally describe precisely the density in $V^1L_M(\Omega)\cap C^{0,\gamma}(\Omega)$, $\gamma\in [0,1]$, is actually elementary. The mentioned trick of~\cite{bgs} is substituted with a different inequality, which is intermediate to the mentioned one. We estimate pointwise the gradient of convolution of our function with a mollifier by a $C^{0,\gamma}$-norm of our function and a parameter of mollification raised to a relevant power, see Lemma~\ref{lem:inqSdinf}. This simple change is surprisingly powerful. Further use of topological arguments makes our procedure independent of any particular structure of the functional, which typically needs to be controlled, cf.~\cite{bacomi-cv, bgs-arma,DeFilippisLeonetti+2022, eslev}. In particular, unlike~\cite[Theorem~4]{bacomi-cv}, we do not make use of any iterative procedure. Consequently, in the case of~\eqref{def-double-phase-functional}, no decay of moduli of the continuity of $b$ with respect to any of its variables is required. In fact, we infer the absence of Lavrentiev's phenomenon for any $\cF$ given by~\eqref{cF-def} with $F:\Omega\times\R\times\rn\to\R$ which is continuous with respect to the last two variables and such that~\eqref{sandwich}, as long as $M$ satisfies one of the prescribed balance conditions. In the case of~\eqref{def-double-phase-functional} it means that $b$ can be taken merely bounded with respect to $x$ and continuous with respect to $u$. Moreover, convexity is imposed only on $\xi\mapsto M(\cdot,\xi)$, not on $\xi\mapsto F(\cdot,\cdot,\xi)$. It is remarkable how simple is the proof of Theorem~\ref{theo:lavrentiev} while having Theorem~\ref{theo:approx}. The mentioned methods allow providing counterparts of these results for vector-valued maps by almost the same arguments, see Theorems~\ref{theo:lavrentiev-vec} and~\ref{theo:approx-vec}. \newline

\noindent{\bf Organization.} 
The framework of Musielak--Orlicz and Musielak--Orlicz--Sobolev spaces, as well as notation and other general information, is exposed in Section~\ref{sec:prelim}. The proofs of the main results, that is Theorems~\ref{theo:lavrentiev} and~\ref{theo:approx}, are given in Section~\ref{sec:proofs}. Section~\ref{sec:general-growth} provides extended results on the absence of Lavrentiev's phenomenon for functionals of general growth and full anisotropy. Section~\ref{sec:vec} is devoted to the corresponding results for problems with vector-valued maps. In Section~\ref{sec:appl-PDEs}, we present a direct application of the theory of existence in the framework of~\cite{C-b}. We end our contribution with Appendix containing easy computations showing the meaning of $\Bcongi$ and $\Bcongo$ in the special cases.

\section{Preliminaries}\label{sec:prelim}

\noindent {\bf Notation. } Throughout the paper, we assume $\Omega\subset\rn$ is a bounded Lipschitz domain, $n\geq 1$. If $U$ is a fixed set, for $\gamma \in (0,1]$ we denote the H\"older seminorm of a function $h$ as
$$
[h]_{0,\gamma; U} := \sup_{x,y \in U, x \not= y} \, 
\frac{|h(x)-h(y)|}{|x-y|^\gamma}.$$
It is well known that the quantity defined above is a seminorm and when $[h]_{0,\gamma;U}<\infty$, we will say that $h$ belongs to the H\"older space $C^{0,\gamma}(U)$. When clear from the context, we will omit the reference to $U$, i.e.: $[h]_{0,\gamma;U}\equiv [h]_{0,\gamma}$. For functions $f, g$, we write that $f(t) = o(g(t))$, if $\lim_{t\to 0}\frac{f(t)}{g(t)}=0$. \newline

\noindent{\bf Convex functions. } Our main reference for inhomogeneous and anisotropic Musielak--Orlicz and Musielak--Orlicz--Sobolev spaces is~\cite{C-b}. In their definitions, the main role is played by a weak $N$-function.\newline

A~function   $M:\Omega\times\Rn\to\rp$ is called a {\it weak $N$-function} if it satisfies the
following conditions:
\begin{enumerate}
\item $ M$ is a Carath\'eodory's function (i.e. measurable with respect to the first variable and continuous with respect to the second one);
\item  $M(x,0) = 0$ and $\xi\mapsto M(x,\xi)$ is a convex function with respect to $\xi$ for a.a. $x \in \Omega$;
\item   $M(x,\xi) = M(x, -\xi)$ for a.a.  $x \in \Omega$ and all $\xi \in \rn$;
\item   there exist two convex functions $m_1,m_2:\rp\to\rp$ such that {$m_1(0)=0$, $m_1(s)> 0$ for $s>0$, 
} 
 and for a.a. $x\in \Omega$ it holds
\begin{equation*}
m_1(|\xi|)\leq M(x, \xi)\leq m_2(|\xi|).
\end{equation*} 
\end{enumerate} 

\noindent A weak $N$-function is called an {\it $N$-function}, if it is super-linear growth, i.e., it additionally satisfies 
\[\lim_{s\to 0}\frac{m_1(s)}{s}=0=\lim_{s\to 0}\frac{m_2(s)}{s}\quad\text{and}\quad\lim_{s\to\infty}\frac{m_1(s)}{s}=\infty=\lim_{s\to\infty}\frac{m_2(s)}{s}\,.\]
As an example of a weak $N$-function that is not an $N$-function we can indicate $M(x,\xi)=\xi$.

 The {\it complementary~function} $M^*$ (called also Legendre's transform and Young's conjugate) to a function  $M:\Omega\times\rn\to\rp$ is defined by \[M^*(x,\eta)=\sup_{\xi\in\rn}(\xi\cdot\eta-M(x,\xi)),\qquad \eta\in\rn,\ x\in\Omega.\] 
 
 We say that a weak $N$-function $M:\Omega\times\rn\to\r$ satisfies $\Delta_2$ condition close to infinity (denoted $M\in\Delta_2$) if there exist  constants $c,C>0$ and nonnegative integrable function $h:\Omega\to\r$ such that for a.e. $x\in\Omega$ it holds 
\begin{equation*}
 M(x,2\xi)\leq cM(x, \xi)+h(x)\qquad\text{for all}\quad \xi\in\rn:\ \ |\xi|>C.
\end{equation*}

\noindent{\bf Function spaces. } The {\it modular} is a functional defined for measurable functions $\xi:\Omega\to\rn$ by the following formula 
\[\vr_{M;\Omega}(\xi):= \int_\Omega M(x,\xi(x))\,dx.\]

Let $M$ be a weak $N$-function. We deal with the three  {\it Musielak--Orlicz classes of functions}: $\mathcal{L}_M,L_M,E_M$. By ${\cal L}_M(\Omega;\rn)$, we denote the generalized Musielak--Orlicz class is the set of all measurable functions $\xi:\Omega\to\rn$ such that $\vr_{M;\Omega}(\xi)<\infty$. Space ${L}_M(\Omega;\rn)$ is the generalized Musielak--Orlicz space, that is the smallest linear space containing ${\cal L}_M(\Omega;\rn)$, equipped with the Luxemburg norm
$||\xi||_{L_M(\Omega;\rn)}:=\inf\left\{\lambda>0:\vr_{M;\Omega}\left(\tfrac 1\lambda \xi\right)\leq 1\right\}.$ Moreover, ${E}_M(\Omega;\rn)$  is the closure in $L_M$-norm of the set of bounded functions.

Directly from the definition, we see that \[{E}_M(\Omega;\rn)\subset {\cal L}_M(\Omega;\rn)\subset { L}_M(\Omega;\rn).\]
The space ${E}_M(\Omega;\rn)$ coincides with the set of all measurable functions for which $\vr_{M;\Omega}(\tfrac{1}{\lambda}\xi)<\infty$ for every $\lambda>0$. Moreover, unlike $L_M$, space $E_M$  is always separable. 

If $M\in\Delta_2$, then
\[{E}_M(\Omega;\rn)= {\cal L}_M(\Omega;\rn)= {L}_M(\Omega;\rn)\]
and $L_M(\Omega;\rn)$ is separable. When both  $M,M^*\in\Delta_2$  then $L_M(\Omega;\rn)$ is reflexive, see~\cite{C-b}.  

 We define the {\it  Musielak--Orlicz--Sobolev space}  $V^{1}L_M(\Omega)$ and $V^{1}E_M(\Omega)$ as follows
\begin{align*} 
V^{1}L_M(\Omega)&:=\big\{f\in W^{1,1}(\Omega):\ \ \nabla f\in L_M(\Omega;\rn)\big\}\,,\\
 V^{1}E_M(\Omega)&:=\big\{f\in W^{1,1}(\Omega):\ \ \nabla f\in E_M(\Omega;\rn)\big\}\,.
\end{align*}where $\nabla$ stands for distributional derivative. The space is considered endowed with the norm
\[
\|f\|_{V^{1}L_M(\Omega)}:=  \|f\|_{L^1(\Omega)}+ \|\nabla f\|_{L_M(\Omega;\rn)}\,.
\]
Zero-trace versions of the spaces will be defined as follows
\begin{align*} 
V^{1}_0L_M(\Omega)&:=\big\{f\in W_0^{1,1}(\Omega):\ \ \nabla f\in L_M(\Omega;\rn)\big\}\,,\\
 V^{1}_0E_M(\Omega)&:=\big\{f\in W_0^{1,1}(\Omega):\ \ \nabla f\in E_M(\Omega;\rn)\big\}\,.
\end{align*}
{If $M\in\Delta_2$, then $L_M(\Omega;\rn)=E_M(\Omega,\rn)$ and we denote $V_0^{1,M}(\Omega):=V^1_0 L_M(\Omega)=V^1_0 E_M(\Omega)$.} \newline

By the Lebesgue's dominated convergence theorem, one can justify the following fact.
\begin{lem}\label{lem:trunc} If $T_k(x)=\min\{k,\max\{-k,x\}\}$ for $k>0$ and $x\in\R$, $M$ is a weak $N$-function, $\vp \in W^{1,1}_0(\Omega)$ and $\nabla \vp\in L_M(\Omega;\rn)$, then for $k\to\infty$ we have $T_k \vp\to\vp$ in $W_0^{1,1}(\Omega)$ and $\nabla T_k\vp\to\nabla \vp$ in $L_M(\Omega;\rn)$.
\end{lem}

\medskip

A typical condition we refer to once studying variable exponent spaces is log-H\"older continuity of the exponent. We say that a function $p$ is log-H\"older continuous, if there exists $c>0$, such that for $x,y$ close enough it holds that
\[|p(x)-p(y)|\leq \frac{c}{\log\left({1}/{|x-y|}\right)}\,.\]
The set of all log-H\"older continuous functions on a bounded set $\Omega$ will be denoted as~$\mathcal{P}^{\log}(\Omega)$. 
\section{Main proofs}\label{sec:proofs}
\subsection{Proof of the density result}

Our approximation in based on the convolution with shrinking.\newline

Let $U$ be a bounded star-shaped domain with respect to a ball $B(x_0, R)$ and let $\kappa_{\delta}^R=1- {\delta}/{R}$. For a~measurable function $\xi:\rn\to\rn$ with $\mathrm{supp}\,\xi\subset U$, we define
\begin{equation}\label{Sdxi}
S_{{U}, \delta}\xi(x) := 
 \int_{U} \rho_\delta( x-y)\xi \left(x_0 + \frac{y - x_0}{\kappa_{\delta}^R}\right)\,dy,
\end{equation}
where $ \rho _\delta(x)=\rho (x/\delta)/\delta^n$ is a standard regularizing kernel on $\rn$  (i.e. $\rho \in C^\infty(\rn)$,
$\mathrm{supp}\,\rho \Subset B(0, 1)$ and $\iO \rho (x)dx = 1$, $\rho (x) = \rho (-x)$, such that $0\leq \rho\leq 1$). Let us notice that $S_{{U},  \delta}(\xi)\in C_c^\infty(\rn;\rn)$. Moreover, for $\delta < \tfrac{R}{4}$ we have $\kappa_{R, \delta}\left( {U} - x_0 \right) + B(x_0, \delta) \subseteq {U}$, which means that for sufficiently small $\delta$ it holds that $S_{{U},  \delta}(\xi) \in C_c^{\infty}({U})$.\newline
Moreover, let us point out that $S_{{U}, \delta}\xi = (\rho_{\delta} * \xi) \left(x_0 + \tfrac{(\cdot) - x_0}{\kappa_{\delta}^R} \right)$, and therefore
\begin{equation}\label{eq:S-grad}
    \nabla S_{{U}, \delta}\xi = \frac{1}{\kappa_{\delta}}S_{{U}, \delta}(\nabla \xi) = (\nabla \rho_{\delta}) * \xi\left(x_0 + \frac{(\cdot) - x_0}{\kappa_{\delta}^R}\right)\,.
\end{equation}
 Let us prove the following lemma.
\begin{lem}\label{lem:S-conv}
Let ${U}$ be a bounded star-shaped domain with respect to a ball $B(x_0, R)$ and let $\xi \in L^1({U})$. Then $S_{{U},  \delta}\,\xi\xrightarrow[\delta \to 0]{} \xi$ in $L^1(U)$ and {in measure}.
\end{lem}
\begin{proof}
Let us firstly abbreviate the notation by setting $S_{\delta} := S_{{U},  \delta}$ and $\kappa_{\delta} := \kappa_{\delta}^R$. Without loss of generality, we assume that ${U}$ is star-shaped with respect to a ball $B(0, R)$. {For a general case one should change variables moving the centre of $B(x_0,R)$ to the origin, then proceed with the proof as below, and then reverse the  change of variables.} {Since convergence in $L^1$ implies convergence in measure} it suffices to show that $\|S_{\delta}\xi - \xi\|_{L^1} \xrightarrow[\delta \to 0]{} 0$. We have
\begin{equation*}
    \|S_{\delta}\xi - \xi\|_{L^1} = \int_{\rn} \left| \int_{B_{\delta}(0)} \rho_{\delta}(y)\xi(\tfrac{x-y}{\kappa_{\delta}})\,dy - \xi(x)  \right|\,dx \leq \int_{\rn} \int_{B_{\delta}(0)} \rho_{\delta}(y)|\xi(\tfrac{x-y}{\kappa_{\delta}})-\xi(x)|\,dy\,dx\,.
\end{equation*}
Let us take a function $g \in C_c^{\infty}(\rn)$. For some ball $B=B(0,r)$ we have
\begin{align*}
\|S_{\delta}g - g\|_{L^1} &\leq \int_{B} \int_{B_{\delta}(0)} \rho_{\delta}(y)|g(\tfrac{x-y}{\kappa_{\delta}}) - g(x)|\,dy\,dx \leq \|\nabla g\|_{L^{\infty}}\int_{B}\int_{B_{\delta}(0)}\rho_{\delta}(y)\left(|x|(\tfrac{1}{\kappa_{\delta}} - 1) + \tfrac{|y|}{\kappa_{\delta}} \right)\,dy\,dx\\
&\leq \|\nabla g\|_{L^{\infty}}\int_{B}\int_{B_{\delta}(0)}\rho_{\delta}(y)\left(r(\tfrac{1}{\kappa_{\delta}} - 1) + \tfrac{\delta}{\kappa_{\delta}} \right)\,dy\,dx \leq \|\nabla g\|_{L^{\infty}}|B|\left(r(\tfrac{1}{\kappa_{\delta}} - 1) + \tfrac{\delta}{\kappa_{\delta}} \right) \xrightarrow[\delta \to 0]{} 0\,.
\end{align*} 
We fix any $\epsilon > 0$ and take $g \in C_c^{\infty}(\rn)$ such that $\|\xi - g\|_{L^1} < \epsilon$. By Young inequality, it holds that
\begin{equation*}
    \|S_{\delta}\xi - \xi\|_{L^1} \leq \|S_{\delta}\xi - S_{\delta}g\|_{L^1} + \|S_{\delta}g-g\|_{L^1} + \|g-\xi\|_{L^1} \leq \kappa_{\delta}^n\epsilon + \|S_{\delta}g-g\|_{L^1} + \epsilon \xrightarrow[\delta \to 0]{} 2\epsilon\,.
\end{equation*}
By taking $\epsilon \to 0$, we obtain $\|S_{\delta}\xi - \xi\|_{L^1} \xrightarrow[\delta\to 0]{} 0$, which implies that $\lim_{\delta\to 0} S_{\delta} \xi=\xi$. %
\end{proof}
Lemma~\ref{lem:S-conv} allows proving the following corollary.
\begin{coro}\label{coro:S-conv}
Let ${U}$ be a bounded star-shaped domain with respect to a ball $B(x_0, R)$ and let $u \in W^{1, 1}_0({U})$. Then $\nabla S_{{U}, \delta}\,\vp \xrightarrow[\delta\to 0]{} \nabla \vp$ in $L^1(U)$ and {in measure}.
\end{coro}
\begin{proof}
Let  $S_{\delta} := S_{{U}, \delta}$ and $\kappa_{\delta} := \kappa_{\delta}^R$. As in Lemma~\ref{lem:S-conv}, without loss of generality, we assume that ${U}$ is star-shaped with respect to a ball $B(0, R)$. Observe that $\nabla S_{\delta}\vp = \tfrac{1}{\kappa_{\delta}}S_{\delta}(\nabla \vp)$ and we have $\|\tfrac{1}{\kappa_{\delta}}S_{\delta}(\nabla \vp) - \nabla \vp\|_{L^1}\leq \|S_{\delta}(\nabla \vp) - \nabla \vp\|_{L^1}+(\tfrac{1}{\kappa_{\delta}} - 1)\|S_{\delta}(\nabla \vp)\|_{L^1}$. Therefore, by Lemma~\ref{lem:S-conv}, it suffices to show that $(\tfrac{1}{\kappa_{\delta}} - 1)\|S_{\delta}(\nabla \vp)\|_{L^1} \xrightarrow[\delta\to 0]{} 0$. 
By Young inequality, we have
\begin{equation*}
    0\leq(\tfrac{1}{\kappa_{\delta}} - 1)\|S_{\delta}(\nabla \vp)\|_{L^1} \leq \kappa_{\delta}^n(\tfrac{1}{\kappa_{\delta}} - 1)\|\nabla \vp\|_{L^1}\|\rho_{\delta}\|_{L^1} \xrightarrow[\delta\to 0]{} 0\,
\end{equation*}
 and the proof is complete. 
\end{proof}
Next, we prove the following inequalities.
\begin{lem}\label{lem:inqSdinf}
Let ${U}$ be a star-shaped domain with respect to a ball $B(x_0, R)$ and $\vp \in W_0^{1, 1}({U})$. It holds that
\begin{itemize}[--]
    \item if $\vp \in L^{\infty}({U})$, then
    \begin{equation}\label{eq:inqLinf}
        \|\nabla S_{U,\delta}(\vp)\|_{L^{\infty}} \leq \delta^{-1}\|\vp\|_{L^{\infty}}\|\nabla \rho\|_{L^1}\,;
    \end{equation}
    \item if $\vp \in C^{0, \gamma}({U})$, $\gamma \in (0, 1]$, then
    \begin{equation}\label{eq:inqH}
        \|\nabla S_{U,\delta}(\vp)\|_{L^{\infty}} \leq \frac{\delta^{\gamma-1}}{\kappa_{\delta}^\gamma}[\vp]_{{0, \gamma}}\|\nabla \rho\|_{L^1}\,.
    \end{equation}
\end{itemize}
\end{lem}
\begin{proof} 
Let  $S_{\delta} := S_{{U},  \delta}$ and $\kappa_{\delta} := \kappa_{\delta}^R$.  As in Lemma~\ref{lem:S-conv} without loss of generality, we assume that $x_0 = 0$. We start with proving~\eqref{eq:inqLinf}. Note that $\nabla S_{\delta}(\vp) = \vp(\tfrac{\cdot}{\kappa_{\delta}}) * (\nabla \rho_{\delta})$. Therefore, by H\"older inequality we have
\begin{equation*}
    \|  \nabla S_{\delta} (\vp) \|_{L^\infty} \leq \|\vp \|_{L^\infty}\int_{\rn}\big|\nabla\rho_\delta(|x|)\big|\,dx = \| \vp \|_{L^\infty}\int_{\rn}\delta^{-n-1}\left|\nabla\rho\left({|x|}/{\delta}\right)\right|\,dx = \delta^{-1}\| \vp \Vert_{L^\infty}\| \nabla \rho\Vert_{L^1}\,,
\end{equation*}
which is the desired result. To prove~\eqref{eq:inqH}, observe firstly that $\nabla S_{\delta}(\vp) = \tfrac{1}{\kappa_{\delta}}S_{\delta}(\nabla \vp)$. Let us fix any $x \in {U}$ and denote $v(y) := \vp(y) - \vp(x/\kappa_{\delta})$. We have
\begin{equation*}
    \nabla S_{\delta}(\vp) = \tfrac{1}{\kappa_{\delta}}S_{\delta}(\nabla \vp) = \tfrac{1}{\kappa_{\delta}}S_{\delta}(\nabla v) = \nabla S_{\delta}(v)\,.
\end{equation*}
Therefore, it holds that
\begin{equation}\label{eq:oct81}
    |\nabla S_{\delta}(\vp)(x)| \leq \int_{{U}} |v(\tfrac{y}{\kappa_{\delta}})||(\nabla \rho_{\delta})(y-x)|\,dy = \int_{B(x, \delta)} |v(\tfrac{y}{\kappa_{\delta}})||(\nabla \rho_{\delta})(y-x)|\,dy \leq \left(\frac{\delta}{\kappa_{\delta}}\right)^{\gamma}[\vp]_{{0, \gamma}}\|\nabla \rho_{\delta}\|_{L^1}\,,
\end{equation}
where we used the fact that $|v(y/\kappa_{\delta})| = |\vp(y/\kappa_{\delta}) - \vp(x/\kappa_{\delta})| \leq \left(\tfrac{|x-y|}{\kappa_{\delta}}\right)^{\gamma}[\vp]_{{0, \gamma}} \leq \left(\frac{\delta}{\kappa_{\delta}}\right)^{\gamma}[\vp]_{{0, \gamma}}$, for $y \in B(x, \delta)$. To end the proof, we observe that $\|\nabla \rho_{\delta}\|_{L^1} = \delta^{-1}\|\nabla \rho\|_{L^1}$ and hence, by~\eqref{eq:oct81} we have the desired result.
\end{proof}
%
%
%
%
%
%
%
%
%
The proof of Theorem~\ref{theo:approx} makes use of the following equivalent modification of $\Bcongi$ as well as $\Bcongo$. 
\begin{description}
    \item {\bf General condition $\Bcong$.} Let us assume that
\[
M(x,\xi)=\sum_{j=1}^L M_j(x,\xi_{j}),
\]
where $M_j:\Omega\times\R^{\ell_j}\to\rp$ are weak $N$-functions, $L>0$, $\xi_j\in\R^{\ell_j}$, and ${\rm lin}_j\{\xi_j\}=\rn$. Suppose there exist constants $c^\Diamond,C^\Diamond>1$ such that 
for every ball $B$ with radius $r \leq 1$, every $x\in B$, all $j$, and for all $\xi_j\in\mathbb{R}^{\ell_j}$ satisfying $(M_j)_{B}^-(C^\Diamond\xi_j) \leq \sup_{|\eta_j|=1}(M_j)_{B}^-(\eta_j c^\Diamond r^{\gamma-1})$, there holds $(M_j)_{B}^+(\xi_j) \leq (M_j)_{B}^-(C^\Diamond\xi_j) + 1\,.$
\end{description}

In this condition, function $M$ can be taken fully anisotropic, even for $L=1$.\newline

We have the following modification of \cite[Theorem~1.2]{ha-aniso}.

\begin{prop}\label{prop:peter}
For a weak $N$-function $M$, the general condition $\Bcong$ implies that there exist constants $C,c>0$ such that 
for every ball $B$ with radius $r \leq 1$, every $x\in B,$ and for all $\xi\in\rn$ it holds \[\MBp(\xi) \leq \Msm(C\xi) + 1\quad\text{ whenever } |\xi| \leq c r^{\gamma - 1}\,.\]
\end{prop}
\begin{proof} First we show that condition $\Bcong$ is equivalent that for every $j$ there exists $c_j,C_j$ such that for every ball $B$ with radius $r \leq 1$, every $x\in B$, all $j$ and for all $\xi_j\in\mathbb{R}^{\ell_j}$  and $(M_j)_{B}^-(C_j\xi_j) \leq \sup_{|\eta_j|=1}(M_j)_{B}^-(\eta_j c_j r^{\gamma-1})$ 
\[\text{there holds}\qquad (M_j)_{B}^+(\xi_j) \leq \big((M_j)_{B}^-\big)^{**}(C_j \xi_j) + 1\,.\] 
The prof of this fact follows almost the same lines as  \cite[Theorem~1.2]{ha-aniso} with the only difference that instead of $\tfrac{K}{|B|}$ we use $\sup_{|\eta_j|=1}(M_j)_{B}^-(\eta_j c_j r^{\gamma-1})$ , that is $s=1+\sup_{|\eta_j|=1}(M_j)_{B}^-(\eta_j c_j r^{\gamma-1})$ . Note that
\begin{equation*}
\frac{s}{s-1} = \frac{\sup_{|\eta_j|=1}(M_j)_{B}^-(\eta_j c_j r^{\gamma-1}) + 1}{\sup_{|\eta_j|=1}(M_j)_{B}^-(\eta_j c_j r^{\gamma-1}) } \leq 1 + \frac{1}{\sup_{|\eta_j|=1}(M_j)_{B}^-(\eta_j c_j r^{\gamma-1}) } \leq 1 + \frac{1}{m_1(c_j)}\,,
\end{equation*}
which means that constants in the proof are bounded independently on $\delta$. {We note that there exists $c$ such that for every $\{\xi_j\colon |\xi_j|\leq c r^{\gamma-1}\}\subseteq\{\xi\colon (M_{j})_B^-(C_j\xi_j) \leq \sup_{|\eta_j|=1}(M_j)_{B}^-(\eta_j c_j r^{\gamma-1})\} $. Therefore, there exists $C$ such that for every $x\in B$
\[M(x,\xi)=\sum_j M_j(x,\xi_j) \leq \sum_j \big((M_j)_{B}^-\big)^{**}(C\xi_j) + 1\leq \Msm(C\xi) +1\quad\text{ whenever }\ {|\xi| \leq c r^{\gamma - 1}}\,.\]
}The last inequality is due to the fact that $\sum_j \big((M_j)_{B}^-\big)^{**}$ is a convex minorant of $M^-_B$.
\end{proof}
\begin{rem}\label{rem:overlineCM}Careful inspection of~\cite{ha-aniso} shows that if condition $\Bcong$ holds for $C^\Diamond$, then $C$ and $c$ from Proposition~\ref{prop:peter} depend only on $C^\Diamond,c^\Diamond,m_1^{-1}(1),m_2^{-1}(1)$.\end{rem}

We are in a position to prove our general result on approximation.
\begin{proof}[Proof of Theorem~\ref{theo:approx}]

We firstly make use of the fact that $\Omega$ is a bounded Lipschitz domain in~$\rn$. That is, by \cite[Lemma~8.2]{C-b}, a set $\overline{\Omega}$ can be covered by a finite family of sets $\{G_i\}_{i=1}^{K}$ such that each $\Omega_i :=\Omega\cap G_i$
is a star-shaped domain with respect to balls $\{ B(x_i, R_i)\}_{i=1}^{K}$, respectively. Then $\Omega=\bigcup_{i=1}^{K}\Omega_i\,.$
Due to~\cite[Proposition 2.3, Chapter 1]{necas}, there exist the partition of unity, i.e., the family $\{\theta_i\}_{i=1}^{K}$ such that
\begin{equation*}
0\le\theta_i\le 1,\quad\theta_i\in C^\infty_c(G_i),\quad \sum_{i=1}^{K}\theta_i(x)=1\ \ \text{for}\ \ x\in\Omega\,.
\end{equation*}
Let us set $R = \min_{i=1}^K R_i$ and $\kappa_{\delta}^{R}$, so that $S_{\Omega_i,\delta}(\theta_i \vp) \in C_c^\infty(\Omega_i)$ for $\delta < \tfrac{R}{4}$. We shall start with proving {\it (i)} simultaneously for $\gamma = 0$ and $\gamma \in (0, 1]$ for bounded functions. For $\gamma = 0$, we take arbitrary $\vp \in V_0^1L_M(\Omega) \cap L^{\infty}(\Omega)$, and for $\gamma \in (0, 1]$, we take $\vp \in V_0^1L_M(\Omega)\cap C^{0, \gamma}(\Omega)$. We define
\begin{equation}\label{eq:Sd-def}
    S_{\delta}(\vp) = \sum_{i=1}^{K} S_{\Omega_i,\delta}(\theta_i \vp)\,.
\end{equation}
Our aim is to show that there exists a constant $\lambda_\vp>0$ such that
\begin{equation}\label{aim}
\lim_{\delta\to 0^+}\iO M\left(x, \frac{\nabla \left(S_\delta(\vp)\right)(x) - \nabla \vp(x)}{\lambda_\vp}\right)\, dx = 0\,.
\end{equation}
By continuity of $M$ and Corollary~\ref{coro:S-conv}, for every $\lambda$ we have  that \[M\left(\cdot, \tfrac{\nabla \left(S_\delta(\vp)\right)(\cdot) - \nabla \vp(\cdot)}{\lambda}\right) \xrightarrow[\delta \to 0]{} 0\quad\text{converge in measure}.\] Therefore, by the Vitali Convergence Theorem, it suffices to show  the family $\left\{M\left(\cdot, \tfrac{\nabla \left(S_\delta(\vp)\right)(\cdot) - \nabla \vp(\cdot)}{\lambda}\right)\right\}_{\delta > 0}$ is uniformly integrable for some $\lambda > 0$. Using the convexity of $M(x, \cdot)$, we obtain
\begin{equation}\label{eq:sep1}
    M\left(x, \frac{\nabla \left(S_\delta(\vp)\right)(x) - \nabla \vp(x)}{\lambda}\right) \leq \frac{1}{2}M\left(x, \frac{2\nabla \left(S_\delta(\vp)\right)(x)}{\lambda}\right) + \frac{1}{2}M\left(x, \frac{2\nabla \vp(x)}{\lambda}\right)\,.
\end{equation}
As $\nabla \vp\in L_M(\Omega;\rn)$, there exists $\lambda_0>0$ such that for $\lambda > \lambda_0$, it holds that $\int_{\Omega} M\left( \frac{2}{\lambda} \nabla \vp(x)\right)\,dx < \infty$. Therefore, by~\eqref{eq:sep1}, it suffices to show that the family $\{M\left(\cdot, \frac{1}{\lambda}\nabla \left(S_\delta(\vp)\right)(\cdot)\right)\}_{\delta > 0}$ is uniformly integrable for some $\lambda > \lambda_0$. 

We recall~\eqref{eq:Sd-def}. Since Jensen's inequality yields
\begin{equation}\label{eq:sep5} M\left(x, \tfrac{1}{\lambda}\nabla \left( S_\delta(\vp) \right)(x)\right) dx \leq \tfrac{1}{K}\sum_{i=1}^K M\left(x, \tfrac{K}{\lambda}\nabla \left(S_{\Omega_i,\delta}(\theta_i \vp)\right)(x)\right)\,, 
\end{equation}
we shall consider integrals $\int_{\Omega_i} M\left(x, \frac{K}{\lambda}\nabla \left(S_{\Omega_i,\delta}(\theta_i \vp)\right)(x)\right)\,dx$ separately.  
    
Now we shall justify that \begin{equation}
        \label{Linftybound}
    \|\nabla S_{\Omega_i,\delta}(\theta_i \vp)\|_{L^{\infty}} \leq C_\vp\delta^{\gamma - 1}
    \end{equation} for some $C_\vp > 0$, which depends on $\vp$ and $\rho$, and for sufficiently small $\delta$. Since we assume that $\vp \in L^{\infty}$, for $\gamma = 0$ we can apply~\eqref{eq:inqLinf} to a function $\theta_i \vp$ to obtain what is needed. For $\gamma \in (0, 1]$, we note that $\theta_i \vp \in C^{0, \gamma}(\Omega)$. This allows us to use~\eqref{eq:inqH} and get that $\|\nabla S_{\Omega_i,\delta}(\theta_i \vp)\|_{L^{\infty}} \leq \frac{\delta^{\gamma-1}}{\kappa_{\delta}^\gamma}[\theta_i \vp]_{{0, \gamma}}\|\nabla \rho\|_{L^1}$. Since for $\delta < \tfrac{R}{4}$ it holds that $\kappa_{\delta} \geq \tfrac{1}{2}$, we have~\eqref{Linftybound}. \newline  
    
    Let $C, c > 0$ be constants from Proposition~\ref{prop:peter} and let $\tau := 2(\tfrac{\text{diam}\Omega}{R}+1)$. There exists $\lambda_1>\lambda_0$ such that for $\lambda>\lambda_1$, we have $\left|{K\nabla \left(S_{\Omega_i,\delta}(\theta_i \vp)\right)(x)} \right| \leq (\tau\delta)^{\gamma-1}c\lambda$. We fix such $\lambda$.     Therefore, we have
    \begin{equation}\label{eq:sep2}
        M\left(x, \frac{K}{\lambda}\nabla \left(S_{\Omega_i,\delta}(\theta_i \vp)\right)(x)\right) \leq \MBxd\left(\frac{KC}{\lambda} \nabla \left(S_{\Omega_i,\delta}(\theta_i \vp)\right)(x)\right) + 1\,. 
    \end{equation}
Let us denote $\lambda_2 := \max\left\{\tfrac{\lambda}{KC},\lambda_0\right\}$. We can now use convexity of $\MBxd$ to obtain 
\begin{align}
\MBxd \left( \frac{KC}{\lambda} \nabla \left(S_{\Omega_i,\delta}(\theta_i \vp)\right)(x)\right)  &\leq \MBxd\left(\tfrac{1}{\lambda_2} \int_{B_{\delta}(0)} \rho _\delta(y)\nabla(\theta_i \vp) \left(x_i + \tfrac{ x-y-x_i}{\kappa_\delta} \right) \,dy\right)\nonumber\\
    &\leq \int_{B(0,\delta)} \rho _\delta(y)\MBxd\left(\tfrac{1}{\lambda_2} \nabla (\theta_i \vp) \big(x_i + \tfrac{ x-y-x_i}{\kappa_\delta}\big)  \right)\,dy\,.\label{eq:oct19}
\end{align}
Note that for $y \in B_{\delta}(0)$ and $\delta < \tfrac{R}{4}$, we have
\begin{equation*}
    \left| x - \left(x_i + \tfrac{x-y-x_i}{\kappa_{\delta}} \right)\right| \leq \tfrac{1 - \kappa_{\delta}}{\kappa_{\delta}}|x-x_i| + \tfrac{|y|}{\kappa_{\delta}} \leq 2\delta\left(\tfrac{\text{diam}\Omega}{R}+1\right) = \tau\delta\,.
\end{equation*}
Therefore, by~\eqref{eq:sep2} and~\eqref{eq:oct19}, it holds that
\begin{align*}
      M\left(x, \frac{K}{\lambda}\nabla \left(S_{\Omega_i,\delta}(\theta_i \vp)\right)(x)\right) &\leq \int_{\rn} \rho_\delta(y)M\left(x_i + \tfrac{x-y-x_i}{\kappa_{\delta}}, \tfrac{1}{\lambda_2} \nabla (\theta_i \vp) \big(x_i + \tfrac{ x-y-x_i}{\kappa_\delta}\big) \right)\,dy +1\\
      &= S_{\Omega_i, \delta}\left(M\left(\cdot,\tfrac{\nabla(\theta_i \vp)(\cdot)}{\lambda_2}\right)\right)(x)+1.
\end{align*}
As $\nabla(\theta_i \vp) = (\nabla \theta_i)\vp + \theta_i \nabla \vp$, we have
\begin{align}\label{eq:sep4}
    M(x, \tfrac{1}{\lambda_2}\nabla(\theta_i \vp)(x)) &\leq \tfrac{1}{2}M(x, \tfrac{2}{\lambda_2}(\nabla\theta_i)(x)\vp(x)) + \tfrac{1}{2}M(x, \tfrac{2}{\lambda_2}\theta_i(x)\nabla \vp(x))\nonumber\\
    &\leq \tfrac{1}{2}m_2\left(\tfrac{2}{\lambda_2}\|\nabla \theta_i\|_{L^{\infty}}\|\vp\|_{L^{\infty}}\right) + \tfrac{1}{2}M(x, \tfrac{2}{\lambda_2}\nabla \vp(x))\,.
\end{align}
Since $\lambda_2 \geq \lambda_0$, the right-hand side of~\eqref{eq:sep4} is integrable. Therefore, by Lemma~\ref{lem:S-conv} we have that \[S_{\Omega_i,\delta}\left(M\left(\cdot,\tfrac{\nabla(\theta_i \vp)(\cdot)}{\lambda_2}\right)\right)\xrightarrow[\delta\to 0]{}M\left(\cdot,\tfrac{\nabla(\theta_i \vp)(\cdot)}{\lambda_2}\right)\quad\text{in $L^1(\Omega_i)$}\,.\]
Collecting~\eqref{eq:sep2} and the above observations,  by the Vitali Convergence Theorem, we get that the family $\left\{M\left(\cdot, \frac{2}{\lambda_2}\nabla \left(S_{\Omega_i,\delta} (\vp)\right)(\cdot)\right)\right\}_{\delta > 0}$, is uniformly integrable. In turn, we get the uniform integrability of the right-hand side of~\eqref{eq:sep5}, which in conjunction with~\eqref{eq:sep1} gives us \eqref{aim} with $\lambda_\vp=\lambda_2/2$. Consequently, we obtain both {\it (i)} and {\it (ii)} with the additional assumption that $\vp \in L^{\infty}(\Omega)$. We note that  Lemma~\ref{lem:trunc} allows avoiding this assumption, which is meaningful for $\gamma = 0$. The modular convergence of gradients  and the Poincar\'e inequality gives the claim {\it (i)}. 
\end{proof}

\subsection{Proof of the absence of the Lavrentiev phenomenon}\label{ssec:proof-abs-Lavr}

Once we are equipped with the result on modular density in our anisotropic space (Theorem~\ref{theo:approx}), we are in a position to concentrate on whether the minimizers of functionals like $\cF$ from~\eqref{cF-def} can be obtained as a~limit of smooth maps.

\begin{proof}[Proof of Theorem~\ref{theo:lavrentiev}] Since $C_c^\infty(\Omega)\subset V_0^1L_M(\Omega)$, it holds that $\inf_{u_0+V^1_0L_M(\Omega)}\cF[u]\leq \inf_{u_0+C_c^\infty(\Omega)}\cF[u]\,.$ Let us concentrate on showing the opposite inequality. For an arbitrary $\varepsilon>0$ there exist $u^\ve\in u_0+V_0^1L_M(\Omega)$ such that \[\cF[u^\ve]\leq\inf_{u_0+V^1_0L_M(\Omega)}\cF[u]+\varepsilon\, .\] Therefore it is enough to show that for any $\varepsilon>0$ it holds
\begin{equation}\label{cFleqcF}
\inf_{u_0+C_c^\infty(\Omega)}\cF[u]\leq \cF[u^\ve]\,.
\end{equation}
Let us fix $\varepsilon>0$ and denote $\overline{u}:=u^\ve-u_0$. By Theorem~\ref{theo:approx}, 
there exists $\{u_k\}_{k\geq 1}\subset C_c^\infty(\Omega)$ such that $u_k\to \overline{u}$ strongly in $L^1(\Omega)$ and in measure, and $\nabla u_k\to\nabla \overline{u}$ modularly in $L_M$. Since $F$ is continuous with respect to the second and the third variable, we infer that \[
\text{$F(x,u_k(x)+u_0(x),\nabla u_k(x) +\nabla u_0(x))\xrightarrow[k\to\infty]{} F(x,u^\ve(x),\nabla u^\ve(x))$ in measure.}
\] Therefore, to get~\eqref{cFleqcF}, {by the Vitali Convergence Theorem}, we only need to show that \begin{equation}
    \label{unif-int}\text{the family $\big\{F(x,u_k(x)+u_0(x),\nabla u_k(x) +\nabla u_0(x))\big\}_{k\geq1}$ {is uniformly integrable}.}\end{equation}
Let us note that by assumption \eqref{sandwich}, convexity of $M$ and by $M\in\Delta_2$ there exists $C>0$ such that \begin{align*}
F(x,u_k(x)+u_0(x),\nabla u_k(x) +\nabla u_0(x)) &\leq LM(x,\nabla u_k(x) +\nabla u_0(x)) +L\\
&\leq C M\left(x, {\nabla u_k(x)}\right) + CM\left(x,{\nabla u_0(x)}\right) +L\;
\end{align*}
for every fixed $k\geq 1$. Note that $\vr_{M;\Omega}(\nabla u_0)<\infty$. Moreover, since $\{\nabla u_k\}_{k\geq 1}$ is modularly convergent and $M\in\Delta_2$, we infer that the family $\{M(x,\nabla u_k(x))\}_{k\geq 1}$ {is uniformly integrable}. Thus, \eqref{unif-int} is justified. In turn, we have~\eqref{cFleqcF} and, consequently, \eqref{no-gap-gen} is proven.

By repeating the same procedure for $u\in V^1_0L_M\cap C^{0,\gamma}(\Omega)$ with the use of Theorem~\ref{theo:approx} {\it (ii)} instead of {\it (i)}, one gets~\eqref{no-gap-hold}. 
\end{proof}

\section{Results for functionals of more general growth}\label{sec:general-growth}

\noindent{\bf Full anisotropy.} Let us recall condition $\Bcong$ formulated right before Proposition~\ref{prop:peter}. We remind that full anisotropy is covered by this condition already with the choice of $L=1$. This is the generality under which the method of the proof of Theorem~\ref{theo:lavrentiev} is directly valid.
\begin{theo} 
\label{theo:lavrentiev-gen}
Let $\Omega$ be a bounded Lipschitz domain in $\rn$ and functional $\cF$ be given by~\eqref{cF-def} with $F:\Omega\times\R\times\rn\to\R$ satisfying~\eqref{sandwich} for a weak $N$-function $M$ {that is continuous with respect to both variables and such that $M\in\Delta_2$}. Assume further that $F$ is measurable with respect to the first variable and continuous with respect to the second and the third variable. Then we observe the absence of Lavrentiev's phenomenon in the following cases.
\begin{enumerate}[{\it (i)}]
    \item If $\gamma=0$, $u_0\in V^1L_M(\Omega)$, and $M$ satisfies condition $\Bcong$, then\begin{align*}
    \inf_{u_0+V^1_0L_M(\Omega)}\cF[u]=\inf_{u_0+C_c^\infty(\Omega)}\cF[u]\,.
\end{align*} 
    \item If $\gamma\in(0,1]$, $u_0\in V^1L_M(\Omega)\cap C^{0,\gamma}(\Omega)$, and $M$ satisfies condition $\Bcong$, then
\begin{align*}
    \inf_{u_0+V^1_0 L_M(\Omega)\cap C^{0,\gamma}(\Omega)}\cF[u]=\inf_{u_0+C_c^\infty(\Omega)}\cF[u]\,.
\end{align*} 
\end{enumerate}
\end{theo} 

\noindent{\bf General growth.} We shall stress an essential difference between the formulation of results, if one does not assume $M\in\Delta_2$. Indeed, this complicates the choice of the natural function space for the problem. Since~\cite{Gossez} it is known that in general the smooth approximation in the classical Orlicz--Sobolev spaces is possible only with respect to the so-called modular topology (not in norm). Under the doubling regime, sequential modular and norm closures coincide. Hence, if $M\in\Delta_2$ we have $L_M=E_M$, $V^1_0 L_M=V^1_0 E_M=V_0^{1,M}$ and the ambiguity disappears. Let us point out that our main results are valid actually without this structure, and that we allow for a natural gap between $V^1_0 L_M$ and $V^1_0 E_M$. In order to relax the growth of $M$ to allow for $M\not\in\Delta_2$, condition $\Bcong$ needs to be strengthened to the following condition. 
\begin{description}
\item {\bf Condition $\Bcongp$.} Let us assume that 
\[
M(x,\xi)=\sum_{j=1}^L M_j(x,\xi_{j}),
\]
where $M_j:\Omega\times\R^{\ell_j}\to\rp$ are weak $N$-functions, $L>0$, $\xi_j\in\R^{\ell_j}$, and ${\rm lin}_j\{\xi_j\}=\rn$. Suppose there exists a constant $C^\Diamond>1$ and a function $\vartheta(r)=\frac{1}{o(r^{1-\gamma})}$, such that 
for every ball $B$ with radius $r \leq 1$, every $x\in B$, all $j$, and for all $\xi_j\in\mathbb{R}^{\ell_j}$ satisfying $(M_j)_{B}^-(C^\Diamond\xi_j) \leq \sup_{|\eta_j|=1}(M_j)_{B}^-\left(\eta_j \vartheta(r)\right)$, there holds $(M_j)_{B}^+(\xi_j) \leq (M_j)_{B}^-(C^\Diamond\xi_j) + 1\,.$
\end{description}
Let us note again, that full anisotropy is covered by this case in particular with the choice of $L=1$.

Then, the non-doubling version of Theorem~\ref{theo:lavrentiev} reads as follows. 
\begin{theo} 
\label{theo:lavrentiev-non-doubling}
Let $\Omega$ be a bounded Lipschitz domain in $\rn$ and functional $\cF$ be given by~\eqref{cF-def} with $F:\Omega\times\R\times\rn\to\R$ satisfying~\eqref{sandwich} for a weak $N$-function $M$ {that is continuous with respect to both variables}. Assume further that $F$ is measurable with respect to the first variable and continuous with respect to the second and the third variable. Then we observe the absence of Lavrentiev's phenomenon in the following cases.
\begin{enumerate}[{\it (i)}]
    \item If $\gamma=0$, $u_0\in V^1E_M(\Omega)$, and $M$ satisfies condition $\Bcongp$, then\begin{align*}
    \inf_{u_0+V^1_0E_M(\Omega)}\cF[u]=\inf_{u_0+C_c^\infty(\Omega)}\cF[u]\,.
\end{align*} 
    \item If $\gamma\in(0,1]$, $u_0\in V^1E_M(\Omega)\cap C^{0,\gamma}(\Omega)$, and $M$ satisfies condition $\Bcongp$, then
\begin{align*}
    \inf_{u_0+V^1_0 E_M(\Omega)\cap C^{0,\gamma}(\Omega)}\cF[u]=\inf_{u_0+C_c^\infty(\Omega)}\cF[u]\,.
\end{align*} 
\end{enumerate}
\end{theo}
\begin{proof}
The proof of Theorem~\ref{theo:lavrentiev-non-doubling} follows the same lines as the proof of Theorem~\ref{theo:lavrentiev} with the only difference that in the lack of $M\in\Delta_2$, it needs to be proven that the approximate sequence $\{\nabla u_n\}$ converge modularly for any $\lambda>0$. This is equivalent to the fact that $\nabla u_n \to \nabla\overline{u}$ strongly in $L_M$. To get this, one needs to use a slightly modified version of Theorem~\ref{theo:approx}. We observe that under $\Bcongp$ for any $\lambda_0$ there exists $\delta$ small enough that \eqref{eq:sep2} holds and 
further modular convergence holds for any $\lambda>0$. In turn, we get that $C_c^\infty$ is dense in the norm topology in $V^1E_M(\Omega)$ provided $\Bcongp$ and the rest of  the reasoning of the proof of Theorem~\ref{theo:lavrentiev} applies under the regime of Theorem~\ref{theo:lavrentiev-non-doubling}.
\end{proof}

\begin{rem}\rm \label{rem:gen-growth}
Note that the example of an essentially fully anisotropic function $\Phi$ from~\cite[Example~1]{cn}  is of the form as in $\Bcongp$. This function $\Phi$ does not satisfy $\Delta_2$, but it is trapped between $|\xi|^p$ and $|\xi|^p\log(1+|\xi|)$. We consider $M(x,\xi)=|\xi|^p+a(x)\Phi(\xi)$. In turn, we can infer that  Lavrentiev's phenomenon for $\cF$ does not occur between $W^{1}_0 E_M (\Omega)$ and $C_c^\infty(\Omega)$ whenever $0\leq a\in {\cal P}^{\log}(\Omega)$. 

We can also consider Orlicz double-phase, that is $M(x, \xi) = \phi(|\xi|) + a(x)\psi(|\xi|)$, without an assumption that Young functions $\phi, \psi \in \Delta_2$. Then, the Lavrentiev's phenomenon for $\mathcal{F}$ does not occur if modulus of continuity of $a$ satisfies $\omega_a(t) \leq \vartheta(o(t^{\gamma - 1}))$, where $\vartheta(t)= \tfrac{\phi(t)}{\psi(t)}$. The proof follows the same line as in~\eqref{eq:orlicz-ex}.
\end{rem}

\begin{rem}\label{rem:pleqn}
\rm
Note that one can obtain results similar to Theorem~\ref{theo:lavrentiev}, \ref{theo:lavrentiev-gen}, and \ref{theo:lavrentiev-non-doubling} on the absence of Lavrentiev's phenomenon with the use of~\cite[Theorem~1]{Bor-Chl} under the regime of $(\mathsf{B})$ therein. We observe that in sub-regions of $\Omega$ where a weak $N$-function $M(x,\xi)$ has a growth not higher than $|\xi|^{\tfrac{n}{1-\gamma}}$, $\gamma\in[0,1]$, there exists $C>0$ such that 
\begin{equation}\label{eq:setinset}
    \left\{\xi\in\rn\colon M_B^{-}(C^\Diamond\xi)\leq \sup_{\eta:\, |\eta|=1}M_B^-(\eta r^{\gamma-1})\right\}\subseteq \left\{\xi\in\rn\colon M_B^{-}(C \xi)\leq r^{-n}\right\}.
    \end{equation}
This implies that condition $\Bcong$ follows from $(\mathsf{B})$ from \cite{Bor-Chl}. Therefore, Theorem~\ref{theo:approx} implies \cite[Theorem~1]{Bor-Chl}.\\
On the other hand, in the sub-regions where $M(x,\xi)\geq |\xi|^{\tfrac{n}{1-\gamma}}$, the opposite inclusion to \eqref{eq:setinset}  is true. Hence, $(\mathsf{B})$ is implied by $\Bcong$. Consequently, one can formulate a mixed condition distinguishing these separate regions. 
\end{rem}
\section{Vector-valued maps}\label{sec:vec}
Let us consider a variational functional
\begin{align*}
    {\boldsymbol\cF}[\bu]:=\int_\Omega F(x,\bu,\nabla \bu)\,dx\,,
\end{align*}
over an open and bounded set $\Omega\subset\rn$, $n,d\geq 1$, where $F:\Omega\times\R^d\times\rnd\to\rp$ is merely continuous with respect to the second and the third variable. Its growth is given by the means of a weak  $N$-function $M:\Omega\times\rnd\to\rp$ that is continuous with respect to both variables and has the following decomposition. We denote $\bxi=(\bxi^1,\dots,\bxi^d)\in\rnd$, $\bxi^i\in\rn$, $i = 1, 2, \dots, d$, and we define a family of weak $N$-functions $M_i : \Omega \times \rn \to [0, \infty)$ such that
\begin{align}\label{M-vec}
    M(x,\bxi):=\sum_{i=1}^d M_i(x,\bxi_i)\,.
\end{align} We suppose that there exist constants $0<\nu,\beta<1<L$ such that
\begin{align*}
    \nu M\left(x,{\beta} \bxi\right)\leq F(x,{\bf z},\bxi)\leq L( M\left(x,\bxi\right)+1)\,, \qquad\text{for all }\ x\in\Omega,\ {\bf z}\in\R^d,\ \bxi\in\rnd\,.
\end{align*}
We recall that $\Bcongi$ or $\Bcongo$ are special cases of $\Bcong$, so the general case we present in particular embraces the case of Theorem~\ref{theo:lavrentiev}. Our most general result reads as follows.

\begin{theo}[Absence of Lavrentiev's phenomenon for vector valued maps]
\label{theo:lavrentiev-vec}
Let $\Omega$ be a bounded Lipschitz domain in $\rn$ and functional ${\boldsymbol\cF}$  for a weak $N$-function $M$ that is continuous with respect to both variables and such that $M\in\Delta_2$ and~\eqref{M-vec} are satisfied.   Then we observe the absence of Lavrentiev's phenomenon in the following cases.
\begin{enumerate}[{\it (i)}]
    \item If  $\gamma=0$, ${\bf u_0}\in V^1L_M(\Omega)$, and $M$ satisfies condition $\Bcong$, then\begin{align*}
    \inf_{{\bf u_0}+V^1_0L_M(\Omega)}{\boldsymbol\cF}[{\bf u}]=\inf_{{\bf u_0}+C_c^\infty(\Omega)}{\boldsymbol\cF}[{\bf u}]\,. 
\end{align*} 
    \item If $\gamma\in(0,1)$, ${\bf u_0}\in V^1L_M(\Omega)\cap C^{0,\gamma}(\Omega)$, and $M$ satisfies condition $\Bcong$, then
\begin{align*}
    \inf_{{\bf u_0}+V^1_0L_M(\Omega)\cap C^{0,\gamma}(\Omega)}{\boldsymbol\cF}[{\bf u}]=\inf_{{\bf u_0}+C_c^\infty(\Omega)}{\boldsymbol\cF}[{\bf u}]\,.
\end{align*} 
    \item If ${\bf u_0}\in V^1L_M(\Omega)\cap C^{0,1}(\Omega)$, then
\begin{align*}
    \inf_{{\bf u_0}+V^1_0L_M(\Omega)\cap C^{0,1}(\Omega)}{\boldsymbol\cF}[{\bf u}]=\inf_{{\bf u_0}+C_c^\infty(\Omega)}{\boldsymbol\cF}[{\bf u}]\,. 
\end{align*} 
\end{enumerate}
\end{theo}
The proof of the above theorem bases on the following approximation result.

\begin{theo}[Density of smooth functions for vector-valued maps]\label{theo:approx-vec}
Let $\Omega$ be a bounded Lipschitz domain in $\rn$ and  $M$ be a weak $N$-function that is continuous with respect to both variables and such that~\eqref{M-vec} is satisfied. Then the following assertions hold true. \begin{enumerate}[{\it (i)}]
\item If $\gamma=0$ and $M$ satisfies condition $\Bcong$, then for any $\bvp\in V_0^1L_M(\Omega)$  there exists a sequence $\{\bvpd\}\subset C_c^\infty(\Omega)$, such that  $ \bvpd\to  \bvp$ strongly in $L^1(\Omega)$ and in measure, and $\nabla\bvpd\xrightarrow[ ]{M}\nabla \bvp$ modularly in $L_M(\Omega;\rnd)$.
\item If  $\gamma\in(0,1]$ and $M$ satisfies condition $\Bcong$, then  for any $\bvp\in V_0^1L_M(\Omega)\cap C^{0,\gamma}(\Omega)$  there exists a sequence $\{\bvpd\}\subset C_c^\infty(\Omega)$, such that  $ \bvpd\to  \bvp$ strongly in $L^1(\Omega)$ and in measure, and $\nabla\bvpd\xrightarrow[ ]{M}\nabla \bvp$ modularly in $L_M(\Omega;\rnd)$.
\end{enumerate} 
Moreover, in both above cases, if $\bvp\in L^\infty(\Omega;\r^d)$, then there exists $c=c(\Omega)>0$, such that $\|\bvpd\|_{L^\infty(\Omega;\r^d)}\leq c \|\bvp \|_{L^\infty(\Omega;\r^d)}$.
\end{theo}

\begin{proof}[On the proofs of Theorems~\ref{theo:lavrentiev-vec} and~\ref{theo:approx-vec}]
Let us take any $\gamma \in [0, 1]$ and consider an $N$-function $M : \Omega \times \rnd \to [0, \infty)$ satisfying $\Bcong$. 
Observe that for every $i$, the function $M_i$ is an $N$-function satisfying condition $\Bcong$. Let us take any $\bvp : \Omega \to \rd$ such that $\bvp \in V^1L_M(\Omega) \cap C^{0, \gamma}(\Omega)$ and let us denote $\bvp := (\vp^1, \vp^2, \dots, \vp^d)$. By Theorem~\ref{theo:approx}, for every $i$ there exists a sequence $\{\vp^{i}_{ \delta}\}_{\delta > 0} \subset C_c^{\infty}(\Omega)$ such that $\vp^i_{ \delta} \xrightarrow[\delta \to 0]{} \vp^i$ strongly in $L^1(\Omega)$ and in measure, and $\nabla \vp^i_{\delta} \xrightarrow[\delta \to 0]{M_i} \nabla \vp^i$ modularly in $L_{M_i}(\Omega;\rn)$. Let us consider a sequence $\bvpd := (\vp^{1}_{\delta}, \vp^{2}_{\delta}, \dots, \vp^{d}_{\delta})$. Note that $\bvpd\xrightarrow[\delta \to 0]{} \bvp$ strongly in $L^1(\Omega)$ and in measure. We shall prove that $\nabla \bvpd \xrightarrow[\delta \to 0]{M} \nabla \bvp$ modularly in $L_M(\Omega; \rnd)$. Let us take any $\lambda > 0$ such that
\begin{equation*}
    \int_{\rn} M_i\left(x, \frac{\nabla \vp_{\delta}^i - \nabla \vp^i}{\lambda/d}\right) \,dx\xrightarrow[\delta \to 0]{} 0, \qquad \text{for $i = 1, 2, \dots, d$.}
\end{equation*}
By~\eqref{M-vec}, we have that
\begin{equation*}
    M\left(x, \frac{\nabla \bvpd - \nabla \bvp}{\lambda} \right)= \sum_{i=1}^{d} M_i\left(x, \frac{\nabla \vp_{ \delta}^i - \nabla \vp_i}{\lambda}\right)\,.
\end{equation*}
Therefore, it holds that 
\begin{equation*}
    \int_{\rnd} M\left(x, \frac{\nabla \bvpd - \nabla \bvp}{\lambda} \right)\,dx \xrightarrow[\delta \to 0]{} 0\,.
\end{equation*}
Hence, we have convergence $\nabla \bvpd \xrightarrow[\delta \to 0]{M} \nabla \bvp$. The absence of Lavrentiev's phenomenon for $\boldsymbol{\cF}$ follows by the same arguments as in the scalar case, see the proof of Theorem~\ref{theo:lavrentiev}.
\end{proof}

\begin{rem}{\rm Note that one can also formulate a counterpart of Theorem~\ref{theo:lavrentiev-non-doubling}  for vector-valued maps.}\end{rem}
\section{Applications to PDEs}\label{sec:appl-PDEs}

One of the direct consequences of our main results is a supplement to the current state of theory of existence of solutions to boundary value problems that are posed in the anisotropic Musielak-Orlicz spaces, which are equipped with modular density of smooth functions. For an illustration, let us consider a second-order elliptic PDE of a form
\[\begin{cases}-\dv \mathcal{A}(x,\nabla u) +B(x,u,\nabla u)=f\quad \text{in }\ \Omega,\\
u=0\qquad \text{on }\ \partial\Omega,\end{cases} \]
where the leading part  of the operator $\mathcal{A}$ is a monotone vector field, being a Caratheodory's function that satisfies growth and coercivity conditions prescribed by the means of an $N$-function $M:\Omega\times\rn\to[0, \infty)$ reading
\begin{align*}
M(x,c_1\xi)&\leq \mathcal{A}(x,\xi)\cdot\xi+h_1(x),\\
c_2 M^*(x,c_3\mathcal{A}(x,\xi))&\leq M(x,c_4\xi)+h_2(x),\end{align*}
for some constants $c_1,c_2,c_3,c_4>0$, fixed $h_1,h_2\in L^1(\Omega)$, almost all $x\in\Omega$ and all $\xi\in\rn{}$. Note that the appearance of all of these constants matters, as there is no kind of doubling condition considered. As it is not important from the point of view of the application of Theorem~\ref{theo:approx}, we skip conditions on $B$ that need to be imposed to obtain existence. Such problems with a regular datum $f$ in the class of weak solutions as well as for merely integrable $f$ in the class of renormalized solutions were studied in~\cite{C-b,CKL,gszg}. Despite the formulations of the theorems there prescribe more restrictive assumptions on an $N$-function $M$, they are used only in order to ensure the modular approximation properties of the space $V^1_0L_M(\Omega)$. In fact, it is enough to make use of Theorem~\ref{theo:approx} {\it (i)} precisely in the form it is stated for example
\begin{itemize}[--]
\item instead of \cite[Theorems~2.2]{gszg} (requiring the condition (M))  to justify the proofs of the results on the existence \cite[Theorems~2.1 and~1.1]{gszg} are proven under condition $\Bcongi$ or $\Bcongo$ for $\gamma=0$;
\item instead of \cite[Theorem 3.7.7]{C-b} (requiring one of the conditions ($Me$) and ($Me$)$_p$) to justify the proofs of the results on the existence  \cite[Theorems~4.1.5 and~5.2.3]{C-b} under condition $\Bcongi$ or $\Bcongo$ for $\gamma=0$.
\end{itemize}
Similarly, one can relax assumptions in~\cite{CKL} and, after minor adaptations, in the case of parabolic problems of \cite{bgs,ren-para,cgzg,Li} and \cite[Chapters~4.2]{C-b}. A bit more challenging is the application of our ideas in the proof of existence of renormalized solutions to parabolic problems when $M$ defining the space is changing not only in space but also along time, i.e., $M=M(t,x,\xi)$, see~\cite{ren-t-para} and \cite[Chapter~5.3]{C-b}. 

\section*{Appendix} 
Let us present the computations for special instances of $M$ satisfying conditions $\Bcongi$ and $\Bcongo$. This not only justifies Corollary~\ref{coro:lavrentiev}, but also indicates the spaces equipped with modular density of smooth functions due to Theorem~\ref{theo:approx}.\newline

\noindent{\bf Examples for $\Bcongi$, $\gamma\in [0,1]$.} Let $M:\Omega\times\rp\to\rp$. We assume additionally that $M\in\Delta_2$.  In all the following cases, it suffices to justify that the balance condition holds for $1\leq |\xi|\leq r^{\gamma-1}$. Indeed, we note that for $\xi$ with sufficiently small $|\xi|$ we have $M(x, \xi) \leq m_2(|\xi|) \leq 1\leq 1+M_B^-(\xi)$. Note that if $|\xi| \in [c, 1]$ for some $c > 0$, then always $M_B^{+}(\xi) \leq m_2(1) \leq m_1(C) \leq m_1(\tfrac{C}{c}|\xi|) \leq M_B^-(\tfrac{C}{c}\xi)$, for $C = m_1^{-1}(m_2(1))$.  Since $M \in \Delta_2$, $\Bcongi$ holds if $\tfrac{M(x, \xi)}{M(y, \xi)} \leq C$ for some $C > 0$, and for all $x, y \in B$ and $|\xi| \in [1, r^{\gamma-1}]$. 

Therefore, let us assume that $1 \leq |\xi| \leq r^{\gamma - 1}$, $|x-y| < r$, and we consider
\begin{enumerate}[{\it (i)}]
    \item $M(x, \xi) = |\xi|^{p(x)}$. Then we have
    \begin{equation*}
        \frac{M(x, \xi)}{M(y, \xi)} = |\xi|^{p(x)-p(y)} \leq r^{(\gamma - 1)|p(x)-p(y)|} \leq C\,, 
    \end{equation*}
     whenever $p \in \mathcal{P}^{\text{log}}$ and $1 \leq p(\cdot) < \infty$.
    \item $M(x, \xi) = |\xi|^p + a(x)|\xi|^p\log(e + |\xi|)$. Then
    \begin{equation*}
        \frac{M(x, \xi)}{M(y, \xi)} = \frac{1 + a(x)\log(e + |\xi|)}{1 + a(y)\log(e+|\xi|)} \leq 1 + |a(x)-a(y)|\log(e + cr^{\gamma-1}) \leq C\,,
    \end{equation*}
    whenever $0 \leq a \in \mathcal{P}^{\text{log}}$ and $1 \leq p<\infty$.
    \item $M(x, \xi) = |\xi|^p + a(x)|\xi|^q$. Then we have
    \begin{equation*}
       \frac{M(x, \xi)}{M(y, \xi)} = \frac{1 + |\xi|^{q-p}a(x)}{1 + |\xi|^{q-p}a(y)} \leq 1 + |\xi|^{q-p}|a(x)-a(y)| \leq 1 + r^{-\alpha}|a(x)-a(y)| \leq C\,,
    \end{equation*}
    whenever $0 \leq a \in C^{0, \alpha}$ and $1 \leq p \leq q \leq p + \tfrac{\alpha}{1 - \gamma}$.
    \item $M(x, \xi) = |\xi|^{p(x)} + a(x)|\xi|^{q(x)}$. Then
    \begin{align*}
        \frac{M(x, \xi)}{M(y, \xi)} = |\xi|^{p(x)-p(y)} \frac{1 + a(x)|\xi|^{q(x)-p(x)}}{1 + a(y)|\xi|^{q(y)-p(y)}} &\leq C\left(1 + |a(x)-a(y)||\xi|^{q(x)-p(x)} + |\xi|^{q(x)-q(y)}|\xi|^{p(x)-p(y)} \right)\\
        &\leq C + C|a(x)-a(y)|r^{-\alpha} \leq C\,,
    \end{align*}
    whenever $0 \leq a \in C^{0, \alpha}$, $p, q  \in \mathcal{P}^{\text{log}}$, and $1 \leq p(\cdot) \leq q(\cdot) \leq p(\cdot) + \tfrac{\alpha}{1 - \gamma}$.
    \item $M(x, \xi) = |\xi|^p + \sum_{i=1}^{n} a_i(x)|\xi|^{q_i}$. Then
    \begin{align*}
        \frac{M(x, \xi)}{M(y, \xi)} = \frac{1 + \sum_{i=1}^{n}a_i(x)|\xi|^{q_i-p}}{1 + \sum_{i=1}^{n}a_i(y)|\xi|^{q_i-p}} \leq \sum_{i=1}^{n} \frac{1 + a_i(x)|\xi|^{q_i - p}}{1 + a_i(y)|\xi|^{q_i - p}} \leq n + \sum_{i=1}^{n}|a_i(x)-a_i(y)|r^{-\alpha_i} \leq C\,,
    \end{align*}
    whenever $0 \leq a_i \in C^{0, \alpha_i}$ and $1 \leq p \leq q_i \leq p + \tfrac{\alpha_i}{1 - \gamma}$, for all $i$.
    \item $M(x, \xi) = \phi(|\xi|) + a(x)\psi(|\xi|)$ for  some $\phi, \psi \in \Delta_2$. Then
    \begin{equation}\label{eq:orlicz-ex}
        \frac{M(x, \xi)}{M(y, \xi)} = \frac{1 + a(x)\frac{\psi(|\xi|)}{\phi(|\xi|)}}{1 + a(y)\frac{\psi(|\xi|)}{\phi(|\xi|)}} \leq 1 + \omega_a(|\xi|^{\frac{1}{\gamma - 1}})\frac{\psi(|\xi|)}{\phi(|\xi|)} \leq C\,,
    \end{equation}
    whenever $\omega_a(t) \leq \tfrac{\phi(t^{\gamma-1})}{\psi(t^{\gamma-1})}$. Note that we comment what happens without assuming $\phi, \psi \in \Delta_2$ in Section~\ref{sec:general-growth}.
    \end{enumerate}

\noindent{\bf Examples for $\Bcongo$, $\gamma\in[0,1]$.}  Let $M(x,\xi)=\sum_{i=1}^n M_i(x,|\xi_{i}|)$,
where $M_i:\Omega\times\rp\to\rp$ are weak $N$-functions and $\xi=(\xi_1,\dots,\xi_n)$. We assume additionally that $M_i\in\Delta_2$, $i=1,\dots,n$. As in the isotropic case, it suffices to justify that the balance condition holds for $1\leq |\xi_i|\leq r^{\gamma-1}$. 
Therefore, let us assume that $1 \leq |\xi_i| \leq r^{\gamma - 1}$, $|x-y| < r$, and for all $i$ we consider
\begin{enumerate}[{\it (i)}]
    \item $M_i(x, |\xi_i|) = |\xi_i|^{p_i(x)}$. Then
    \begin{equation*}
        \frac{M_i(x, |\xi_i|)}{M_i(y, |\xi_i|)} \leq |\xi_i|^{p_i(x)-p_i(y)} \leq C_i\,,
    \end{equation*}
    whenever $1 \leq p_i \in \mathcal{P}^{\text{log}}$.
    \item $M_i(x, |\xi_i|) = |\xi_i|^{p_i} + a_i(x)|\xi_i|^{q_i}$. Then
    \begin{equation*}
        \frac{M_i(x, |\xi_i|)}{M_i(y, |\xi_i|)} = \frac{1 + a_i(x)|\xi_i|^{q_i-p_i}}{1 + a_i(y)|\xi_i|^{q_i-p_i}} \leq 1 + |a_i(x)-a_i(y)||\xi_i|^{q_i-p_i} \leq 1 + |a_i(x)-a_i(y)|r^{-\alpha_i} \leq C_i\,,
    \end{equation*}
    whenever $0 \leq a_i \in C^{0, \alpha_i}$ and $1 \leq p_i \leq q_i \leq p_i + \tfrac{\alpha_i}{1 - \gamma}$.
    \item $M_i(x, |\xi_i|) = |\xi_i|^{p_i} + a_i(x)|\xi_i|^{p_i}\log(e + |\xi_i|)$. Then
    \begin{equation*}
        \frac{M_i(x, |\xi_i|)}{M_i(y, |\xi_i|)} \leq \frac{1 + a_i(x)\log(e + |\xi_i|)}{1 + a_i(y)\log(e + |\xi_i|)} \leq 1 + |a_i(x)-a_i(y)|\log(e + cr^{\gamma - 1}) \leq C_i\,,
    \end{equation*}
    whenever $0 \leq a_i \in \mathcal{P}^{\text{log}}$.
    \item $M_i(x, |\xi_i|) = |\xi_i|^{p_i(x)} + a_i(x)|\xi_i|^{q_i(x)}$. Then
    \begin{align*}
        \frac{M_i(x, |\xi_i|)}{M_i(y, |\xi_i|)} &\leq |\xi_i|^{p_i(x)-p_i(y)} \frac{1 + a_i(x)|\xi_i|^{q_i(x)-p_i(x)}}{1 + a_i(y)|\xi_i|^{q_i(y)-p_i(y)}}\\
        &\leq C\left(1 + |a_i(x)-a_i(y)||\xi_i|^{q_i(x)-p_i(x)} + |\xi_i|^{q_i(x)-q_i(y)}|\xi_i|^{p_i(x)-p_i(y)} \right)\\
        &\leq C\left(1 + |a_i(x)-a_i(y)|r^{-\alpha_i}\right) \leq C_i\,,
    \end{align*}
    whenever $a_i \in C^{0, \alpha_i}$, $p_i, q_i \in \mathcal{P}^{\text{log}}$, and $1 \leq p_i(\cdot) \leq q_i(\cdot) \leq p_i(\cdot) + \tfrac{\alpha_i}{1 - \gamma}$.
\end{enumerate}

\section*{Acknowledgement}
  The authors express gratitude to Peter H\"ast\"o (University of Turku) for a discussion on the meaning of balance conditions during the Thematic Research Programme Anisotropic and Inhomogeneous Phenomena in Warsaw in September 2022. The authors would like to thank Piotr Rybka (University of Warsaw) and Cristiana De Filippis (University of Parma) for insighting comments.
  
{\small  
\bibliography{bib}
} 
\end{document}